\documentclass[12pt,reqno]{amsart}
\usepackage{amsmath,amsthm,amssymb,amsfonts,amscd}
\usepackage{mathrsfs}
\usepackage{bbm}
\usepackage{bbding}

\usepackage{bm}
\usepackage{hyperref}

\usepackage{geometry}
\usepackage{color}
\usepackage{xcolor}

\usepackage{picture,epic}
\usepackage{tikz}

\numberwithin{equation}{section}

\theoremstyle{plain}
\newtheorem{theorem}{Theorem}[section]
\newtheorem{lemma}[theorem]{Lemma}
\newtheorem{corollary}[theorem]{Corollary}

\theoremstyle{definition}

\theoremstyle{remark}
\newtheorem{remark}[theorem]{Remark}

\renewcommand{\Re}{\operatorname{Re}}
\renewcommand{\Im}{\operatorname{Im}}
\newcommand{\vol}{\operatorname{vol}}

\newcommand{\sgn}{\operatorname{sgn}}

\newcommand{\GL}{\operatorname{GL}}
\newcommand{\SL}{\operatorname{SL}}

\renewcommand{\mod}{\operatorname{mod}\ }

\makeatletter
\def\@tocline#1#2#3#4#5#6#7{\relax
  \ifnum #1>\c@tocdepth 
  \else
    \par \addpenalty\@secpenalty\addvspace{#2}%
    \begingroup \hyphenpenalty\@M
    \@ifempty{#4}{%
      \@tempdima\csname r@tocindent\number#1\endcsname\relax
    }{%
      \@tempdima#4\relax
    }%
    \parindent\z@ \leftskip#3\relax \advance\leftskip\@tempdima\relax
    \rightskip\@pnumwidth plus4em \parfillskip-\@pnumwidth
    #5\leavevmode\hskip-\@tempdima
      \ifcase #1
       \or\or \hskip 1em \or \hskip 2em \else \hskip 3em \fi%
      #6\nobreak\relax
    \hfill\hbox to\@pnumwidth{\@tocpagenum{#7}}\par
    \nobreak
    \endgroup
  \fi}
\makeatother

\begin{document}

\title[Uniform subconvexity bounds for $\GL(2)\times \GL(2)$ in the spectral aspect]
{Uniform subconvexity bounds for $\GL(2)\times \GL(2)$
\\ $L$-functions in the spectral aspect}

\author{Zhao Xu}

\address{School of Mathematics \\ Shandong University \\ Jinan \\ Shandong 250100 \\China}
\email{zxu@sdu.edu.cn}


\date{\today}

\begin{abstract}
  In this paper, we study the second moment for $\GL(2)\times \GL(2)$ $L$-functions
  $L(\frac{1}{2},f\times g)$,
  which leads to a uniform subconvexity bound in the spectral aspect.
In particular, if either $f$ or $g$ is a dihedral Maass newform, or if one of them has level
$1$, we obtain a Burgess-type bound that is uniform in both
$t_f$ and $t_g$, where $t_f$, $t_g$ denote the spectral parameters of $f$, $g$.
As an application, we also establish a shrinking result for QUE in the case of dihedral Maass newforms.
\end{abstract}



\thanks{}

\maketitle

\section{Introduction} \label{sec:Intr}

One of the important problems in the theory of $L$-functions is to ask for
good upper bounds on the critical line which are stronger than
the ones obtained by the Phragmen--Lindel$\ddot{\textup{o}}$f principle.
Starting with the work from Weyl \cite{Weyl21}, many considerable results
have been achieved (see \cite{Bur63}, \cite{Sar01}, \cite{Mic04},
\cite{JutMot05}, \cite{HarMic06},
\cite{MicVen10}, \cite{Li11}, \cite{Mun15III},
\cite{Mun15IV}, \cite{PetYou20}, \cite{BloJanNel}, \cite{HuMicNel}, etc).
In particular, the subconvexity problem of the Rankin--Selberg $L$-functions
and the triple $L$-functions has received much attention.
In this paper, we consider the subconvexity bound for $GL(2)\times GL(2)$ $L$-function
uniformly in the spectral aspect.

Let $f$ and $g$ be two holomorphic or Maass cusp newforms.
Denote by $t_f$ the weight or the spectral parameter of $f$.
Let $L(s,f\times g)$ be the associated Rankin-Selberg $L$-function.
The first subconvexity result of $L(s,f\times g)$
was obtained by Sarnak \cite{Sar01} (in the weight aspect),
who proved that (for $f$ and $g$ holomorphic),
\begin{align*}
  L\left(\frac{1}{2}+it,f\times g\right)\ll
   t_f^{\frac{18}{19-2\theta}+\varepsilon},
\end{align*}
where $\theta$ is any exponent towards
the Ramanujan--Petersson conjecture, and $\theta\leq 7/64$
due to Kim--Sarnak \cite{Kim2003}.
The subconvex exponent was improved by Lau--Liu--Ye \cite{LauLiuYe06}
to the Weyl-type result $2/3$, where the levels of $f$ and $g$
are coprime or equal, and the nebentypus are trivial.
Both results of Sarnak and Lau--Liu--Ye are based on the spectral analysis
and estimates for triple products of automorphic forms.

Before Lau--Liu--Ye's work, Blomer \cite{Blomer05} obtained
\begin{align*}
  L\left(\frac{1}{2}+it,f\times g\right)\ll_{g,D,t,\varepsilon}
   t_f^{\frac{6-2\theta}{7-4\theta}+\varepsilon},
\end{align*}
for $f$ and $g$ of any level and any nebentypus.
Blomer's approach is via
Jutila's variant of the circle method to detect the shifted condition
and the Kuznetsov trace formula to exploit cancellation for
sums of Kloosterman sums.

The first hybrid subconvexity for $L(s,f\times g)$
was obtained by Jutila--Motohashi \cite{JutMot06} (in the $t_f$- and $t$-aspects),
who obtained that, if $D=q=1$,
\begin{align*}
  L\left(\frac{1}{2}+it,f\times g\right)
  \ll_{g,\varepsilon}
  \begin{cases}
    t_f^{\frac{2}{3}+\varepsilon}, & 0\leq t\ll t_f^{\frac{2}{3}},
    \\
    t_f^{\frac{1}{2}}t^{\frac{1}{4}}, & t_f^{\frac{2}{3}}\leq t\ll t_f,
    \\
    t^{\frac{3}{4}+\varepsilon}, & t_f\ll t\ll t_f^{\frac{3}{2}-\varepsilon}.
  \end{cases}
\end{align*}
It is reasonable that their method is also effective for the
Hecke congruence subgroups.

Furthermore, Michel--Venkatesh \cite{MicVen10} obtained the uniform subconvexity
bound for $L(1/2,f\times g)$ in all aspects of $f$.
In a recent preprint, Nelson \cite{Nel21} obatined an immense breakthrough,
which addresses the spectral aspect
for all standard $L$-functions on $GL_n$ in the case of uniform parameter growth.

We now turn to our notation.
Let $\chi$ and $\psi$ be two Dirichlet characters modulo $D$ and $q$, respectively,
where $D$ and $q$ are two positive integers.
Let $\mathcal{B}^*(D,\chi)$ be the set of $L^2$-normalised Maass newforms
on $\Gamma_0(D)$ with nebentypus $\chi$ and spectral parameters $t_f$,
where $\Gamma_0(D)$ is the Hecke congruence subgroup of $SL_2(\mathbb{Z})$.
For $g\in \mathcal{B}^*(q,\psi)$, the $L^4$-norm of $g$ is defined by
\begin{align*}
  ||g||_4=\left(\int_{\Gamma_0(q)\backslash\mathbb{H}}
  |g(z)|^4d\mu z\right)^{\frac{1}{4}},
\end{align*}
where $\mathbb{H}$ is the upper half plane and $d\mu z=\frac{dxdy}{y^2}$.
Let $T$ be a large parameter, and denote by $T_0=|T+t_g|$ and $L=|T-t_g|$.

\begin{theorem}\label{thm:moment result}
Let notation be as above. Assume $q$ is square free and $\psi$ is a real primitive character
modulo $q$.
Let $H$ be a parameter such that $T_0\ll (TH)^{\frac{3}{4}-\varepsilon}$
and $H\ll \min\{T^{1-\varepsilon},L/\log T\}$.
Then, we have
\begin{align}\label{eqn:main moment 1}
  \sum_{\substack{f\in\mathcal{B}^*(D,\chi)\\ T-H\leq t_f\leq T+H}}
  \left|L\left(\frac{1}{2},f\times g\right)\right|^2\ll
  T_0^{1+\varepsilon}L^{\frac{1}{2}}||g||_4^2+T^{1+\varepsilon}H.
\end{align}
Moreover, if $t_g\ll \frac{T^{1+\varepsilon}}{H}$, we have
\begin{align}\label{eqn:main moment 2}
  \sum_{\substack{f\in\mathcal{B}^*(D,\chi)\\ T-H\leq t_f\leq T+H}}
  \left|L\left(\frac{1}{2},f\times g\right)\right|^2\ll
  T^{\frac{3}{2}+\varepsilon}+T^{1+\varepsilon}H.
\end{align}
\end{theorem}

  Note that the moment results in Theorem \ref{thm:moment result} are not
  Lindel$\ddot{\textup{o}}$f hypothesis on average results.
  However, they are enough for our purpose (Corollary \ref{cor:subconvex}).
  Furthermore, it is reasonable to assume $t_g$ is real and larger than
  a big constant, since otherwise
  Lau--Liu--Ye
  \cite{LauLiuYe06} has obtained a stronger moment result which implies
  the Weyl-type subconvexity bound.
  In Theorem \ref{thm:moment result}, we need $\psi$ to be real due to
  technical reasons (see \eqref{eqn:nead to be real} below).

  We turn to the relation among our parameters.
  Firstly, it is natural to assume $H\leq T^{1-\varepsilon}$.
  Secondly, the condition $T_0\ll (TH)^{\frac{3}{4}-\varepsilon}$ comes from
  the uniform asympotic formula for the Bessel function $J_{2it}(x)$,
  which we will use twice
  (after the  Kuznetsov trace formula
  and after the Voronoi summation formula).
  We also need this condition when dealing with the hypergeometry function
  $F$, where we use it to make sure the Dirichlet series expression \eqref{eqn:def F}
  is absolutely and uniformly convergent in our crucial range
  (see \S \ref{subsec:deal hypergeometry}).
  Hence, in our results, we have the basic set up
  $t_g\leq T_0\ll (TH)^{\frac{3}{4}-\varepsilon}\ll T^{\frac{3}{2}-\varepsilon}$
  and $H\gg T_0^{\frac{4}{3}+\varepsilon}T^{-1}\gg T^{\frac{1}{3}+\varepsilon}$.
  It is worth mentioning that the lower bound of $H$ ($\gg T^{\frac{1}{3}+\varepsilon}$ in the
  special case of $t_g\ll T$) also occured in the work of Lau--Liu--Ye \cite{LauLiuYe06}
  and in the very recent Weyl bound work for the triple product $L$-functions of
  Blomer--Jana--Nelson \cite{BloJanNel}.
  Lastly, we need $H=o(L)$ to keep the analytic conductor of all the $L$-functions $L(\frac{1}{2},f\times g)$
  in Theorem \ref{thm:moment result} being of size $T_0^2L^2$.

We will use the method as in Sarnak \cite{Sar01} and Lau--Liu--Ye \cite{LauLiuYe06}.
A key feature in \cite{LauLiuYe06} is to provide the analytic continuation
of the shifted convolution sum to a bigger region and finally use
the deep results in \cite{Goo81,Jut96,Jut97,KroSta04} (\cite[(3.14)]{LauLiuYe06})
to treat certain triple products.
Our basic observation is that, rather than applying \cite[(3.14)]{LauLiuYe06},
one can use the recent $L^4$-norm
to control the triple products of automorphic forms.
Actually, by the Watson--Ichino formula, the triple products of automorphic
forms can be reduced to moment of triple product $L$-functions,
by which one can start to consider the $L^4$-norm problem.
It turns out that the $L^4$-norm conjecture ($||g||_4\ll t_g^\varepsilon$)
would yield the Burgess bound uniformly in the spectral aspect.

Furthermore, the Bessel functions coming from $g$ require a careful treatment.
We apply the uniform expansion of the
$J$-Bessel function, a standard tool in several previous works
(see \cite{BloJanNel}, \cite{BloBut19}, \cite{BloBut20}).


The hypergeometric function
$F$ also appears naturally in our process and needs a precise evaluation.
When $g$ is fixed,
$F$ can be replaced by $1$ with an acceptable error (see \cite[\S 12]{LauLiuYe06}),
but since $g$ is varying here, we provide a more detailed analysis.
The main approach is outlined in \S \ref{subsec:sketch of the main thm},
with full details given in \S \ref{subsec:deal hypergeometry}-\S \ref{subsec:treat P2}.

Now, we combine the $L^4$-norm results together with Theorem \ref{thm:moment result}
to get our subconvexity bounds.

  For $g$ being a dihedral Maass newform, Luo \cite{Luo} got
  \begin{align*}
    ||g||_4\ll t_g^{\varepsilon},
  \end{align*}
  when $q$ is prime,
  while Humphries--Khan \cite{HumKhan20} derived the asymptotic formula
  \begin{align}\label{eqn:HK l4 d}
    ||g||_4^4=\frac{3}{\textup{vol}(\Gamma_0(q)\backslash \mathbb{H})}
    +O_q(t_g^{-\delta}),
  \end{align}
  where $\delta>0$ is an absolute positive constant and $q$ is square free.

\begin{corollary}\label{cor:subconvex}
  Let conditions be as in Theorem \ref{thm:moment result}.
  Furthermore, let $D$ be square free and $\chi$ be a real primitive character modulo $D$.
  Then, if one of $f$ and $g$ is a dihedral Maass newform, we have the
  Burgess-type bound
    \begin{align}\label{eqn:subconvex burgess}
    L\left(\frac{1}{2},f\times g\right)\ll |t_f+t_g|^{\frac{3}{4}+\varepsilon}.
  \end{align}
\end{corollary}
Recently, Ki \cite{Ki} announced that
\begin{align}\label{eqn:Ki}
    ||g||_4\ll t_g^{\varepsilon},
  \end{align}
  when $q=1$.
For $g\in \mathcal{B}^*(q)$
\footnote{We omit $\psi$ in $\mathcal{B}(q,\psi)$ and $\mathcal{B}^*(q)$
if $\psi$ is trivial.}, where $q$ is square free,
Humphries--Khan \cite{HumKhan22} proved that
  \begin{align}\label{eqn:HK l4}
    ||g||_4\ll t_g^{\frac{3}{152}+\varepsilon}.
  \end{align}

  In Theorem \ref{thm:moment result}, the assumption that $\psi$ is primitive can be dropped.
We impose primitivity only to apply the explicit version of the Watson-Ichino formula given by Humphries–-Khan \cite{HumKhan20}; this formula in fact holds in general setting (with level-dependence left implicit).
Hence we obtain the following result.
  \begin{corollary}\label{cor:subconvex2}
  Let $f\in \mathcal{B}^*(D)$ and $g\in \mathcal{B}^*(q)$. Then, we have
  \begin{align}\label{eqn:subconvex level=1}
   L\left(\frac{1}{2},f\times g\right)\ll |t_f+t_g|^{\frac{3}{4}+\frac{3}{152}+\varepsilon}.
 \end{align}
 In particular, if one of $D$ and $q$ is $1$, we have the Burgess-type bound
 \begin{align}\label{eqn:subconvex burgess level=1}
    L\left(\frac{1}{2},f\times g\right)\ll |t_f+t_g|^{\frac{3}{4}+\varepsilon}.
  \end{align}
  \end{corollary}

\begin{remark}
  Let notation be as in Corollary \ref{cor:subconvex} or \ref{cor:subconvex2}.
  We can also obtain a subconvexity from \eqref{eqn:main moment 1}
  when $T_0^{\frac{1}{3}+\delta}\leq L=o(T_0) $ for any positive number $\delta$.
  Clearly, we have $t_f\asymp T\asymp t_g\asymp T_0$ now.
  Let conditions be as in Corollary \ref{cor:subconvex}.
  Note that the convexity bound of $L(\frac{1}{2},f\times g)$ is
  $T^{\frac{1}{2}+\varepsilon}L^{\frac{1}{2}}$.
  Then, by using \eqref{eqn:HK l4}
  and \eqref{eqn:HK l4 d},
  one has $T_0^{1+\varepsilon}L^{\frac{1}{2}}||g||_4^2=o(T^{1+\varepsilon}L)$.
  So \eqref{eqn:main moment 1} gives a subconvexity bound by letting
   $H=LT^{-\frac{\delta}{2}}$.
\end{remark}

\subsection{Sketch of the proof of Theorem \ref{thm:moment result}}
\label{subsec:sketch of the main thm}
The overall outline follows \cite{LauLiuYe06}.
By the approximate functional equation in \cite{LiYou12}
(see Lemma \ref{lemma:AFE-LiYoung}),
the estimation in Theorem \ref{thm:moment result} is reduced to
\begin{align*}
  \sum_{\substack{f\in B^*(D,\chi)\\ T-H\leq t_f\leq T+H}}
  \left|\sum_{n\geq 1}\lambda_f(n)\lambda_g(n)V\left(\frac{n}{N}\right)\right|^2,
\end{align*}
where $\lambda_f(n)$ and $\lambda_g(n)$ denote the $n$-th
Hecke eigenvalues of $f$ and $g$, $V$ denotes a nice function
(see \S \ref{subsec:smooth weight functions}),
and $N\ll T_0^{1+\varepsilon}L$ does not vary with $f$.
After the Kuznetsov trace formula and the Voronoi summation formula,
one needs to treat the $J$-Bessel functions $J_{2it}(2x)$ and $J_{2it_g}(2u)$.
We follow the technique in \cite{BloJanNel}:
apply the uniform expansion of the $J$-Bessel function (see \eqref{eqn:J2it})
and Lemma \ref{lemma:stationary_phase BKY} to truncate at $x$, $u\gg T^{1-\varepsilon}H$,
and use the Taylor expansion \eqref{eqn:Taylor of omega} to give
an explicit expression of $J$.
Then we need to deal with an integral of the shape (see \eqref{x integral})
\begin{align*}
  \int_{\mathbb{R}}V\left(\frac{x}{N}\right)
  e(\alpha x^{1/2}+\beta x^{-1/2})dx,
\end{align*}
which has been clearly computed in \cite{BloJanNel}
by using Lemma \ref{lem:main theorem in KPY}
(the stationary phase method).
So one sees that, after these treatment, we arrive at
\begin{align*}
  \sum_{r\asymp R}P(r,N)
  \sum_{\substack{c\asymp C\\ c\equiv0(\mod[q,D])}}
  G_{\chi\psi}(r,c),
\end{align*}
where $C\ll \frac{N}{T^{1-\varepsilon}H}$ and $R=\frac{T_0LC^2}{N}$
($\approx C^2$ in the ``generic" case).
Here $G_{\chi\psi}(r,c)$ is a character sum and
$P(r,N)$ is a shifted convolution sum.
By the inverse Mellin transform, we have
\begin{align*}
  P(r,N)=\frac{1}{4\pi i}\int_{(2)}(2N)^s\tilde{G}_r(s)D_{g}(s,1,1,r)ds,
\end{align*}
where $D_g$ is the shifted convolution sum \eqref{eqn:scs}
(with $\nu_1=\nu_2=1$)
and $\tilde{G}_r(s)$, in any large but fixed vertical stripe,
is negligible unless $|\Im s|\asymp V_0$ with $V_0=\frac{T_0LC}{N}$
($\approx C$ in the ``generic" case),
in which case $\tilde{G}_r(s)$ can be replaced by
an explicit expression (see \eqref{eqn: replace G}).

For $D_{g}(s,1,1,g)$, we need to relate it to $D_{g,F}(s,1,1,g)$ (see \eqref{eqn:scs+F})
which has an extra hypergeometric $F$,
so that we can use the spectral decomposition
(see \eqref{eqn:decomposition of D}-\eqref{eqn:DE}).
Our situation is different from \cite{Sar01} and \cite{LauLiuYe06},
since $g$ is varying.
In fact, we can replace $D_g$ with $D^\dag_g$ here (see \eqref{eqn:express P by D}),
where ``$\dag$" means that the $m$ in $D_g$ or $D_{g,F}$
is restricted to $|m|\gg NT^{-\varepsilon}$.
In our crucial ranges (see \eqref{recall parameters}), one can approximate
$D^\dag_{g}(s,1,1,g)$ by a linear combination of $D_{g,F}^{\dag}(s+2\ell,1,1,r)$
(see \eqref{eqn:express D dag}), where $\ell \leq A$ with $A$ being an absolutely
large constant.
One can replace $D^\dag_{g,F}$ by $D_{g,F}$, if the contribution of
$D^\ddag_{g,F}:=D_{g,F}-D^\dag_{g,F}$ is acceptable.
In fact, by some identities of special functions together with
integration by parts and shifting the contour to far left,
we can show that the contribution $D^\ddag_{g,F}$ can be omitted.

To handle $D_{g,F}$, we first use the spectral decomposition.
Then, after truncating $t_\phi$ and $\tau$ in \eqref{eqn:Dd} and \eqref{eqn:DE},
respectively, we arrive at
\begin{multline*}
  \sum_{1\leq k\ll \log T}\sum_{r\asymp R}\mathcal{S}_{r,k}
  \sum_{\substack{\phi\in \mathcal{B}^*(q_1)
  \\ t_\phi\asymp V_0}}
  \frac{\overline{\rho}_\phi(r)}{\cosh(\frac{\pi t_\phi}{2})}\left<\phi,|g|^2\right>_q
  \int_{\substack{\Re s=2\ell+2\\|\Im s|\asymp V_0}}
  N^s\frac{\tilde{G}_r(s-2\ell)}{r^{s-1}}
  \left\{\cdots\right\}ds
  \\
  +\textup{similar\ part\ of\ the\ continuous\ spectrum},
\end{multline*}
where $q_1|q$, $\mathcal{S}_{r,k}$ ($\ll C^{1+\varepsilon}$) is a certain sum,
and the essential ingredients in $\{\cdots\}$
are the gamma factors in the spectral decomposition
of $D_{g,F}$ (see \eqref{eqn:Dd}).
We move the $s$-integral to $\Re s=-\frac{1}{2}$
passing possible poles.
For simplicity, we only explain the idea of dealing with the contribution
coming from the discrete spectrum of these poles,
since all the other terms can be treated similarly.
Now we use \eqref{eqn: replace G} to replace $\tilde{G}_r$,
and see that $r$ and $\Im s$ in the factor
$\frac{\tilde{G}_r(s-2\ell)}{r^{s-1}}$
have been separated completely.
Consequently, by using the Cauchy--Schwarz inequality, we need to estimate
\begin{align*}
  \mathfrak{C}_1
  =\sum_{t_\phi\asymp V_0}\frac{1}{\cosh \pi t_\phi}
  \left|\sum_{r\asymp R}\bar{\rho}_\phi(r)\mathcal{S}_{r,k}\right|^2,
\end{align*}
and
\begin{align*}
  \mathfrak{C}_2=\sum_{t_\phi\asymp V_0}
  e^{\Omega(t_\phi,t_g)}(1+|2t_g-t_\phi|)^{\frac{1}{2}}
  |\left<\phi,|g|^2\right>_q|^2.
\end{align*}
where $\Omega(t_\phi,t_g)$ is defined as in \eqref{eqn:define omega}.
We can get the upper bound for $\mathfrak{C}_1$ through
the large sieve inequality (see Lemma \ref{lem:large sieve}).

For $\mathfrak{C}_2$, if $t_g\gg \frac{T^{1+\varepsilon}}{H}$
(which implies that $t_g\gg V_0^{1+\varepsilon}$ and $\Omega(t_\phi,t_g)=0$),
then, by the spectral decomposition,
it can be controlled by $T_0^{1/2}||g||_4^4$.

If $t_g\ll \frac{T^{1+\varepsilon}}{H}$ (which implies that $L\asymp T_0\asymp T$),
by the Waston--–Ichino formula
(see Lemma \ref{lem:Watson-Ichino q}),
$\mathfrak{C}_2$ can be reduced to treating
\begin{align*}
  \sum_{t_\phi\asymp V_0}L\left(\frac{1}{2},\phi\right)
  L\left(\frac{1}{2},\phi\times \textup{ad}\,g\right).
\end{align*}
In this case, we use \cite[Proposition 6.1]{HumKhan22} to treat the above
mixed $L$-functions, which, roughly speaking, is to perform a dyadic-subdivision,
Cauchy--Schwarz inequality,
the approximate functional equation and the large sieve inequality.
For our purpose, we state the moment result of this special case
($t_g\ll \frac{T^{1+\varepsilon}}{H}$)
as Theorem \ref{thm:moment result}-\eqref{eqn:main moment 2}.



\subsection{Quantum unique ergodicity in shrinking sets}

The quantum unique ergodicity conjecture (QUE) of Rudinick--Sarnak \cite{RudSar94}
states that for if $g(z)$ is a $L^2$-normalized Hecke-Maass newform
and $h$ is a fixed, smooth and compactly-supported function
on $\Gamma_0(q)\backslash \mathbb{H}$, then
\begin{align*}
  \int_{\Gamma_0(q)\backslash \mathbb{H}}
  h(z)|g(z)|^2 d\mu z\rightarrow
  \frac{1}{\textup{vol}(\Gamma_0(q)\backslash \mathbb{H})}
  \int_{\Gamma_0(q)\backslash \mathbb{H}}
  h(z) d\mu z,
\end{align*}
as the Laplace eigenvalue $\lambda_g=\frac{1}{4}+t_g^2$ tends to infinity.
This has been proved via the work of Lindenstrauss \cite{Lin06}
and Soundararajan \cite{Sou10}.
The analog of the QUE conjecture for holomorphic forms (the Mass equidistribution
conjecture) was proved by Holowinsky--Soundararajan \cite{HolSou10}.

A natural question of the QUE conjecture is to determine whether equidistribution
still holds if $\phi(z)$ is supported in a thin set in terms of $t_g$.
To this end, denote by $B=B_R(w)$ be the hyperbolic ball of radius $R$
centred at $w\in \Gamma_0(q)\backslash \mathbb{H}$ with volume
$4\pi \sinh^2(\frac{R}{2})$.
Under the assumption of the generalised Lindel$\ddot{\textup{o}}$f hypothesis,
Young \cite[Proposition 1.5]{You16} obtained that, for $q=1$ and $R\gg t_g^{-\delta}$
with $0<\delta<\frac{1}{3}$, one has
\begin{align}\label{eqn:shrinking statement}
  \frac{1}{\vol(B_R)}\int_{B_R}
  |g(z)|^2 d\mu z\sim
  \frac{1}{\vol(\Gamma_0(q)\backslash \mathbb{H})}.
\end{align}
Moreover, Young also obtained an analogous result for the Eisenstein series
$E(z,\frac{1}{2}+it)$ which states that (\cite[Theorem 1.4]{You16})
\begin{align}\label{eqn:Young shrinking result}
  \frac{1}{\log (1/4+t_g^2)\vol(B_R)}\int_{B_R}
  \left|E\left(z,\frac{1}{2}+it\right)\right|^2 d\mu z\sim
  \frac{1}{\vol(\Gamma_0(q)\backslash \mathbb{H})},
\end{align}
holds unconditionally when $0<\delta<\frac{1}{9}$.
The exponent $\delta$ for \eqref{eqn:Young shrinking result}
has been subsequently improved to $0<\delta<\frac{1}{6}$
by Humphries \cite[Theorem 1.16]{Hum18}.
When $g$ is a a dihedral Maass newform,
Humphries--Khan \cite[Theorem 1.7]{HumKhan20} obatined
an average version which states that
QUE holds for almost every shrinking ball whose radius is larger than
$t_g^{-1}$.

Motivated by these work,
as an application
of Corollary \ref{cor:subconvex},
we prove the following result.
\begin{theorem}\label{thm:shrinkring result for dihedral}
  Let $q\equiv 1(\bmod 4)$ be a positive squarefree fundamental discriminant
  and let $\psi$ be the primitive quadratic character modulo $q$.
  If $R\gg t_g^{-\delta}$ with $0<\delta<\frac{1}{12}$,
  then \eqref{eqn:shrinking statement} holds as $t_g$ tends to infinity along
  any subsequence of dihedral Maass newforms $g\in \mathcal{B}^*(q,\psi)$.
\end{theorem}

$\mathbf{Notation.}$
Throughout the paper, $\varepsilon$ is an arbitrarily small positive number;
$A$ is an sufficiently large positive number.
All of them may be different at each occurrence.
We will also borrow the concepts (``flat" and ``nice") and techneque
in \cite[\S 2.4]{BloJanNel}, from which one can separate the variables in
a nice function (see the details in \S \ref{subsec:smooth weight functions}).



\section{Background}\label{sec:background}

\subsection{Automorphic forms and $L$-functions}\label{subsec:forms and L functions}

We start this subsection by stating some basic concepts from the theory of Maass forms
of weight zero in the context of the Hecke congruence group $\Gamma_0(q)$, where $q$
is a positive integer.
The Petersson inner product is defined by
\begin{align*}
  \left<\phi_1,\phi_2\right>_q=\int_{\Gamma_0(q)\backslash \mathbb{H}}
  \phi_1(z)\overline{\phi_2(z)}d\mu z.
\end{align*}
Let $\mathcal{B}(q,\psi)$ denote an orthonormal basis of Maass cusp forms
of level $q$, nebentypus $\psi$,
with $t_\phi$ denoting the spectral parameter of $\phi\in \mathcal{B}(q,\chi)$.
We have the Fourier expansion
\begin{align}\label{eqn:Fourier expansion}
  \phi(z)=\sum_{n\neq 0}\rho_\phi(n)W_{0,it_\phi}(4\pi|n|y)e(nx),
\end{align}
where $W_{0,it}(y)=(\frac{y}{\pi})^{\frac{1}{2}}K_{it}(\frac{y}{2})$ is a Whittaker function.
If $\phi\in \mathcal{B}^*(q,\psi)$, where $\mathcal{B}^*(q,\psi)\subset \mathcal{B}(q,\psi)$ is the set of $L^2$-normalized newforms
$\phi$ of level $q$, nebentypus $\psi$, then we denote its normalized Hecke eigenvalues
$\lambda_{\phi}(n)$ and record the relation
\begin{align}\label{eqn:Fc1}
  \lambda_{\phi}(n)\rho_{\phi}(1)=\sqrt{n}\rho_{\phi}(n)
\end{align}
for $n\geq 1$, and
\begin{align}\label{eqn:Fc2}
  \rho_{\phi}(-n)=\varepsilon_{\phi} \rho_{\phi}(n),
\end{align}
where $\varepsilon_{\phi}=\pm 1$ is the parity of $\phi$.
Moreover, by Rankin-Selberg theory
and works of Iwaniec \cite{Iwa90} and Hoffstein-Lockhart \cite{HL94},
we have the well-known bounds
\begin{align}\label{eqn:Fc3}
  q^{1-\varepsilon}(1+|t_\phi|)^{-\varepsilon}\ll
  \frac{|\rho_{\phi}(1)|^2}{\cosh(\pi t_{\phi})}
  \ll q^{1+\varepsilon}(1+|t_\phi|)^\varepsilon.
\end{align}
For each singular cusp $\mathfrak{a}$ of $\Gamma_0(q)$ for $\psi$, we define the Eisenstein series
\begin{align*}
  E_{\mathfrak{a}}(z,s,\psi) =
  \sum_{\gamma\in\Gamma_{\mathfrak{a}}\setminus\Gamma_0(q)}
  \bar{\psi}(\gamma)\Im(\sigma_{\mathfrak{a}}^{-1}\gamma z)^s,
\end{align*}
which converges absolutely for $\Re s>1$ and
by analytic continuation for all $s\in\mathbb{C}$, where
$\Gamma_{\mathfrak{a}}$ is the stability group of $\mathfrak{a}$ and
the scaling matrix
$\sigma_{\mathfrak{a}}\in \mathrm{SL}_2(\mathbb{R})$ is such that
$\sigma_{\mathfrak{a}}\infty=\mathfrak{a}$ and
$\sigma_{\mathfrak{a}}^{-1}\Gamma_{\mathfrak{a}}\sigma_{\mathfrak{a}}=\Gamma_\infty$.
The Eisenstein series $E_{\mathfrak{a}}(z,s,\psi)$ is independent of the choice of
scaling matrix, and the Fourier expansion at $s=\frac{1}{2}+i\tau$ can be written as
\begin{align*}
  E_{\mathfrak{a}}\left(z,\frac{1}{2}+i\tau,\psi\right)
  =\delta_{\mathfrak{a}\infty}y^{\frac{1}{2}+i\tau}
  +\varphi_{\mathfrak{a}}\left(\frac{1}{2}+i\tau\right)y^{\frac{1}{2}-i\tau}
  +\sum_{n\neq 0}\rho_{\mathfrak{a}}(n,\tau)W_{0,i\tau}(4\pi|n|y)e(nx).
\end{align*}

To treat the Fourier coefficients, we state the spectral large sieve inequalities
of Deshouillers--Iwaniec \cite{DesIwa}, where the nebentypus $\psi$ is trivial
which is enough for our needs.
Recently, Drappeau \cite{Dra} and  Zacharias \cite{Zac} have extended these results to general $\psi$.
\begin{lemma}\label{lem:large sieve}
  Let $T$, $M\geq 1$, $q\in\mathbb{N}$ and let $(a_m)$, $M\leq m\leq 2M$, be a sequence of
  complex numbers. Then the quantities
  \begin{align*}
    \sum_{\substack{\phi\in \mathcal{B}(q,\psi)\\|t_\phi|\leq T}}\frac{1}{\cosh (\pi t_\phi)}
    \left|\sum_ma_m\sqrt{m}\rho_{\phi}(\pm m)\right|^2,
  \end{align*}
  \begin{align*}
    \sum_{\mathfrak{a}\, \mathrm{singular}}
    \int_{-T}^T\frac{1}{\cosh (\pi t)}\left|\sum_ma_m\sqrt{m}\rho_{\mathfrak{a}}(\pm m,\tau)\right|^2 d\tau
  \end{align*}
  are bounded by
  \begin{align*}
    M^\varepsilon\left(T^2+\frac{M}{q}\right)\sum_m|a_m|^2.
  \end{align*}
\end{lemma}

Now we turn to $L$-functions. Most of the facts can be found in
\cite[Chapter 5]{IwaKow04}.
We write $L(s,\pi)$ for a general $L$-function, which has the Euler product
\begin{align*}
  L(s,\pi)=\prod_p L_p(s,\pi).
\end{align*}
Let $\Lambda(s,\pi):=q(\pi)^{s/2}L_\infty(s,\pi)L(s,\pi)$
denote the completed $L$-function,
where $q_\pi$ denotes the arithmetic conductor of $\pi$,
and $L_\infty(s,\pi)$ is the archimedean part of
$\Lambda(s,\pi)$,
which is of the form $\pi^{-ns/2}\prod_{j=1}^n\Gamma(\frac{s+\kappa_{\pi,j}}{2})$
for some $\kappa_{\pi,j}\in\mathbb{C}$.
The analytic conductor $q(s,\pi)$ of $L(s,\pi)$ is defined as
$q(s,\pi)=q(\pi)\prod_{j=1}^d(|s+\kappa_j|+3)$.
We denote by
\begin{align}\label{eqn:Lq}
L_q(s,\pi):=\prod_{p|q}L_p(s,\pi),\quad L^q(s,\pi):=\frac{L(s,\pi)}{L_q(s,\pi)},\quad
\Lambda^q(s,\pi):=\frac{\Lambda(s,\pi)}{L_q(s,\pi)}.
\end{align}

Let $f\in \mathcal{B}^*(D,\chi)$ and $g\in \mathcal{B}^*(q,\psi)$.
For $\Re s>1$, $L(s,f\times g)$ is given by (see \cite[p. 117]{Blomer05})
\begin{align*}
  L(s,f\times g)=L(2s,\chi\psi)\sum_{n\geq 1}\frac{\lambda_f(n)\lambda_g(n)}{n^s}
  \sum_{d|(qD)^\infty}\frac{\gamma_{f\times g}(d)}{d^s},
\end{align*}
where $\gamma_{f\times g}(d)$ are certain coefficients satisfying
\begin{align*}
  \gamma_{f\times g}(d)\ll_\varepsilon d^{2\theta+\varepsilon}.
\end{align*}
Here $\theta$ stands for an approximation towards the Ramanujan-Petersson conjecture.
The current best result is due to Kim--Sarnak \cite{Kim2003} that $\theta=\frac{7}{64}$
holds.
The Rankin-Selberg $L$-function $L(s,f\times g)$ has analytic continuation to $\mathbb{C}$
and satisfies thd functional equation (see \cite[p. 133]{IwaKow04}
 and \cite[p. 609]{HarMic06})
\begin{align*}
  \Lambda(s,f\times g)=q(f\times g)^{\frac{s}{2}}L_{\infty}(s,f\times g)L(s,f\times g)
  =\varepsilon_{f\times g}\overline{\Lambda(1-\bar{s},f\times g)},
\end{align*}
where $|\varepsilon_{f\times g}|=1$
and
\begin{multline*}
  L_{\infty}(s,f\times g)
  =\pi^{-2s}\Gamma\left(\frac{s+it_f+it_g+\eta}{2}\right)
  \Gamma\left(\frac{s+it_f-it_g+\eta}{2}\right)
  \\
  \cdot\Gamma\left(\frac{s-it_f+it_g+\eta}{2}\right)
  \Gamma\left(\frac{s-it_f-it_g+\eta}{2}\right).
\end{multline*}
Here $\eta=0$, $1$ according to whether $\varepsilon_f=\varepsilon_g$ or not.

We need the following approximate functional equation in Li--Young \cite[Lemma 2.4]{LiYou12}.
\begin{lemma}\label{lemma:AFE-LiYoung}
  Suppose that $q(\frac{1}{2},f\times g)\ll Q$
  for some number $Q>1$. Then there exists a function $W(x)$ depending on $Q$ and $\varepsilon$
  only, such that $W(x)$ is supported on $x\in [0,Q^{\frac{1}{2}+\varepsilon}]$ and satisfies
  \begin{align*}
    x^jW^{(j)}(x)\ll_{j,\varepsilon}1,
  \end{align*}
  where the implied constant depends on $j$ and $\varepsilon$ only $($not on $Q$$)$, such that
  \begin{align}\label{eqn:app}
    \left|L\left(\frac{1}{2},f\times g\right)\right|^2
    \ll Q^\varepsilon \int_{-\log Q}^{\log Q}\left|\sum_{n\geq 1}\frac{\lambda_f(n)\lambda_g(n)}{n^{\frac{1}{2}+iw}}W(n)\right|^2 dw
    +O(Q^{-100}),
  \end{align}
  where the implied constant depends on $\varepsilon$ and $W$ only.
\end{lemma}
\begin{remark}
  Note that $\lambda_f(n)\lambda_g(n)\neq \lambda_{f\times g}(n)$,
  so \eqref{eqn:app} is slightly different from \cite[Lemma 2.4-(2.29)]{LiYou12}.
  But this is unimportant, since one can absorb the factor $L(2s,\chi\psi)\sum_{d|(qD)^\infty}\frac{\gamma_{f\times g}(d)}{d^s}$
  into the weight functions $V_{f,t}$ and $V^*_{f,-t}$ in \cite[Lemma 2.4]{LiYou12},
  and do the same treatment.
\end{remark}


\subsection{Spectral decomposition
of the shifted convolution sum in Sarnak \cite{Sar01}}\label{subsec:shifted sum}

For our purpose, we need the following decomposition result.
\begin{lemma}\label{lem:orthonormal basis}
   $\mathrm{(}$\cite[Lemma 3.1]{HumKhan20}
   $\mathrm{or}$
  \cite[Proposition 2.6]{IwaLuoSar00}$\mathrm{)}$
  An orthonormal basis of the space of Maass cusp forms of
  squarefree level $q$, and trivial nebentypus is given by
\begin{align*}
  \mathcal{B}(q)=\{\phi_\ell: \phi\in\mathcal{B}^*(q_1),
  q_1q_2=q,\ell|q_2\},
\end{align*}
where each newform $\phi\in\mathcal{B}^*(q_1)$ is normalised such that
$\left<\phi,\phi\right>_q=1$ and
\begin{align*}
  \phi_\ell:=\left(L_\ell(1,\textup{sym}^2\phi)\frac{\varphi(\ell)}{\ell}\right)^{\frac{1}{2}}
  \sum_{w_1w_2=\ell}\frac{\nu(w_1)}{w_1}\frac{\mu(w_2)\lambda_\phi(w_2)}{\sqrt{w_2}}\iota_{w_1}\phi.
\end{align*}
Here $\nu(n)=n\prod_{p|n}(1+p^{-1})$, $\varphi(n)=n\prod_{p|n}(1-p^{-1})$,
and $(\iota_{w_1}\phi)(z)=\phi(w_1 z)$.
\end{lemma}

Now we recall the shifted convolution sum in \cite{Sar01}.
Since we are considering the uniform bound in $t_f$ and $t_g$, one needs
a more careful normalization of $g$. Hence, we state the details here.
Let $g(z)\in \mathcal{B}^*(q,\psi)$,
so it has the same Fourier expansion and properties \eqref{eqn:Fourier expansion}-\eqref{eqn:Fc3}
with replaced $\phi$ by $g$.
Let $\nu_1$, $\nu_2$ and $r$ be positive integers.
Denote by
\begin{align}\label{eqn:scs}
  D_{g}(s,\nu_1,\nu_2;r)=\sum_{\substack{m,n\neq0\\\nu_1n-\nu_2m=r}}
  \frac{\overline{\rho_g(m)}\rho_g(n)|mn|^{\frac{1}{2}}}{|\rho_g(1)|^2}
  \left(\frac{\sqrt{\nu_1\nu_2|mn|}}{\nu_1|m|+\nu_2|n|}\right)^{2it_g}
  (\nu_1|m|+\nu_2|n|)^{-s},
\end{align}
and
\begin{multline}\label{eqn:scs+F}
\begin{split}
  D_{g,F}(s,\nu_1,\nu_2;r)=\sum_{\substack{m,n\neq0\\\nu_1n-\nu_2m=r}}
  \frac{\overline{\rho_g(m)}\rho_g(n)|mn|^{\frac{1}{2}}}{|\rho_g(1)|^2}
  \left(\frac{\sqrt{\nu_1\nu_2|mn|}}{\nu_1|m|+\nu_2|n|}\right)^{2it_g}
  (\nu_1|m|+\nu_2|n|)^{-s}
  \\
  \cdot F\left(\frac{s}{2}+it_g,\frac{1}{2}+it_g,\frac{s+1}{2};
  \left(\frac{|\nu_1m|-|\nu_2n|}{|\nu_1m|+|\nu_2n|}\right)^2\right),
\end{split}
\end{multline}
where $F$ is the hypergeometric function.

Let $G(z)=g(\nu_1z)\overline{g(\nu_2z)}$,
and
\begin{align*}
U_r(z,s)=\sum_{\gamma\in\Gamma_\infty\backslash\Gamma}\Im(\gamma z)^se(r\Re(\gamma z)).
\end{align*}
Then, by the standard unfolding method, one has
\begin{multline*}
  \left<U_r(\cdot,s),G\right>_{q\nu_1\nu_2}
  =4(\nu_1\nu_2)^{\frac{1}{2}}\sum_{\substack{m,n\neq 0\\\nu_1m-\nu_2n=r}}
  \overline{\rho(m)}\rho_g(n)|mn|^{\frac{1}{2}}
  \\
  \cdot \int_0^\infty
  y^{s-1}K_{it_g}(2\pi |\nu_1m|y)K_{it_g}(2\pi|\nu_2n|y)dy.
\end{multline*}
The $y$-integral above is (see \cite[6.576-4]{GraRyz})
\begin{multline*}
  \int_0^\infty y^{s-1}K_{it_g}(2\pi|\nu_1m|y)K_{it_g}(2\pi|\nu_2n|y)\textup{d}y
  \\
  =\frac{1}{8\pi^s}\frac{|\nu_1m|^{-it_g-s}|\nu_2n|^{it_g}}{\Gamma(s)}\Gamma\left(\frac{s+2it_g}{2}\right)
  \Gamma\left(\frac{s-2it_g}{2}\right)\Gamma\left(\frac{s}{2}\right)^2
  F\left(\frac{s+2it_g}{2},\frac{s}{2},s;1-\left(\frac{\nu_2n}{\nu_1m}\right)^2\right).
\end{multline*}
Using the transformation for $F$ (see \cite[\S 2.1.5-(24)]{EMOT})
\begin{align}\label{eqn:transform F}
  F\left(a,b,2b;\frac{4z}{(1+z)^2}\right)=(1+z)^{2a}
  F\left(a,a-b+\frac{1}{2},b+\frac{1}{2};z^2\right),
\end{align}
we have
\begin{multline}\label{eqn:eqn of F}
  F\left(\frac{s+2it_g}{2},\frac{s}{2},s;1-\left(\frac{\nu_1n}{\nu_2m}\right)^2\right)
  \\
  =\left(\frac{|\nu_1m|+|\nu_2n|}{2|\nu_1m|}\right)^{-s-2it_g}
  F\left(\frac{s}{2}+it_g,\frac{1}{2}+it_g,\frac{s+1}{2};
  \left(\frac{|\nu_1m|-|\nu_2n|}{|\nu_1m|+|\nu_2n|}\right)^2\right).
\end{multline}
Therefore, we get
\begin{align*}
  \frac{2^{1-s-2it_g} \pi^s (\nu_1\nu_2)^{-\frac{1}{2}}\Gamma(s)}{|\rho_g(1)|^2\Gamma\left(\frac{s+2it_g}{2}\right)
  \Gamma\left(\frac{s-2it_g}{2}\right)\Gamma\left(\frac{s}{2}\right)^2}
  \left<U_r(\cdot,s),G\right>_{q\nu_1\nu_2}=D_{g,F}(s,\nu_1,\nu_2;r).
\end{align*}
On the other hand, by the spectral decomposition, we have
\begin{multline*}
  \left<U_r(\cdot,s), G\right>_{q\nu_1\nu_2}
  =\sum_{\phi\in \mathcal{B}(q\nu_1\nu_2)}
  \left<U_r(\cdot,s), \phi\right>_{q\nu_1\nu_2}
  \left<\phi,G\right>_{q\nu_1\nu_2}
  \\
  +\frac{1}{4\pi}\sum_{\mathfrak{a}}\int_{-\infty}^\infty
  \left<U_r(\cdot,s), E\left(\cdot,\frac{1}{2}+i\tau\right)\right>_{q\nu_1\nu_2}
  \left<E\left(\cdot,\frac{1}{2}+i\tau\right),G\right>_{q\nu_1\nu_2}d\tau.
\end{multline*}
By the unfolding method and \cite[6.561-16]{GraRyz}, one has
\begin{align*}
  \left<U_r(\cdot,s),\phi\right>_{q\nu_1\nu_2}
  =\frac{\pi^{\frac{1}{2}-s}\overline{\rho_{\phi}(r)}}{2r^{s-1}}
  \Gamma\left(\frac{s-\frac{1}{2}+it_\phi}{2}\right)
  \Gamma\left(\frac{s-\frac{1}{2}-it_\phi}{2}\right),
\end{align*}
and
\begin{align*}
  \left<U_r(\cdot,s),E\left(\cdot,\frac{1}{2}+i\tau\right)\right>_{q\nu_1\nu_2}
  =\frac{\pi^{\frac{1}{2}-s}\overline{\rho_{\mathfrak{a}}(r,t)}}{2r^{s-1}}
  \Gamma\left(\frac{s-\frac{1}{2}+i\tau}{2}\right)
  \Gamma\left(\frac{s-\frac{1}{2}-i\tau}{2}\right).
\end{align*}
Therefore, by letting
\begin{align}\label{eqn:Ag}
  A_{g}(s)=\frac{\pi^{\frac{1}{2}}\Gamma(s)}{2^{s+2it_g}|\rho_g(1)|^2\Gamma\left(\frac{s+2it_g}{2}\right)
  \Gamma\left(\frac{s-2it_g}{2}\right)\Gamma^2\left(\frac{s}{2}\right)},
\end{align}
and
\begin{align}\label{eqn:B}
  B(s,\mu)=\Gamma\left(\frac{s-\frac{1}{2}+i\mu}{2}\right)
  \Gamma\left(\frac{s-\frac{1}{2}-i\mu}{2}\right)\cosh\left(\frac{\pi \mu}{2}\right),
\end{align}
we get
\begin{align}\label{eqn:decomposition of D}
  D_{g,F}(s,\nu_1,\nu_2;r)=D_{g,F,d}(s,\nu_1,\nu_2;r)+D_{g,F,E}(s,\nu_1,\nu_2;r),
\end{align}
where
\begin{align}\label{eqn:Dd}
  D_{g,F,d}(s,\nu_1,\nu_2;r)=
  \sum_{\phi\in \mathcal{B}(q\nu_1\nu_2)}
  \frac{\overline{\rho_\phi(r)}(\nu_1\nu_2)^{-\frac{1}{2}}}{r^{s-1}\cosh(\frac{\pi t_\phi}{2})}
  A_{g}(s)B(s,t_\phi)\left<\phi,G\right>_{q\nu_1\nu_2},
\end{align}
and
\begin{multline}\label{eqn:DE}
  D_{g,F,E}(s,\nu_1,\nu_2;r)
  \\
  =\frac{1}{4\pi}\sum_{\mathfrak{a}}\int_{-\infty}^\infty
  \frac{\overline{\rho_{\mathfrak{a}}(r,\tau)}(\nu_1\nu_2)^{-\frac{1}{2}}}
  {r^{s-1}\cosh(\frac{\pi \tau}{2})}
  A_{g}(s)B(s,\tau)
  \left<E_{\mathfrak{a}}\left(\cdot,\frac{1}{2}+i\tau\right),G\right>_{q\nu_1\nu_2}d\tau.
\end{multline}

\subsection{The Watson--Ichino Formula}

In the present work, we only need ``$\nu_1=\nu_2=1$" in the shifted convolution sum $D_{g,F}$,
so we require the Watson--Ichino formula (see \cite{Wat02} and \cite{Ich08}), which relates
$\left<|g|^2,\phi\right>_q$ and $\left<|g|^2,E_\infty(\cdot,\frac{1}{2}+i\tau)\right>_q$
to a triple product $L$-function.
Here we quote the results stated in \cite{HumKhan20}.
\begin{lemma}\label{lem:inner product exchange}
  $\mathrm{(}$\cite[p. 53]{HumKhan20}$\mathrm{)}$
  For $q$ squarefree and $g\in \mathcal{B}^*(q,\psi)$ with $\psi$ primitive, we have
  \begin{align*}
    \sum_{\substack{\phi\in \mathcal{B}(q)\\ t_\phi=\mu}}\left|\left<|g|^2,\phi\right>_q\right|^2
    =\sum_{q_1q_2=q}2^{\omega(q_2)}\frac{\nu(q_2)\varphi(q_2)}{q_2^2}
    \sum_{\substack{\phi\in \mathcal{B}^*(q_1)\\ t_\phi=\mu}}
    \frac{L_{q_2}(1,\mathrm{sym}^2\phi)}{L_{q_2}(\frac{1}{2},\phi)}
    \left|\left<|g|^2,\phi\right>_q\right|^2,
  \end{align*}
  for any $\mu\in [0,\infty)\cup i(0,\frac{1}{2})$.
  Similarly, we have
  \begin{align*}
    \sum_{\mathfrak{a}}
    \left|\left<|g|^2,E_{\mathfrak{a}}\left(\cdot,\frac{1}{2}+i\tau\right)\right>_q\right|^2
    =2^{\omega(q)}
  \left|\left<|g|^2,E_{\infty}\left(\cdot,\frac{1}{2}+i\tau\right)\right>_q\right|^2,
  \end{align*}
  where $\omega(q)=\sharp\{p|n\}$, $\nu(n)=n\prod_{p|n}(1+p^{-1})$,
  and $\varphi(n)=n\prod_{p|n}(1-p^{-1})$.
\end{lemma}
\begin{lemma}\label{lem:Watson-Ichino q}
  $\mathrm{(}$\cite[Proposition 1.16]{HumKhan20}$\mathrm{)}$
  Let $q=q_1q_2$ be squarefree and let $\psi$ be a primitive Dirichlet character modulo $q$.
  Then for $g\in\mathcal{B}^*(q,\psi)$ and for $\phi\in\mathcal{B}^*(q_1)$ of parity
  $\epsilon_\phi\in\{1,-1\}$ normalised such that
  $\left<g,g\right>_q=\left<\phi,\phi\right>_q=1$,
  \begin{align}\label{eqn:Watson-Ichino cusp q}
    |\left<|g|^2,\phi\right>_q|^2=\frac{1+\epsilon_\phi}{16\sqrt{q_1}\nu(q_2)}
    \frac{\Lambda(\frac{1}{2},\phi)\Lambda(\frac{1}{2},\phi\times \textup{ad}\,g)}{\Lambda(1,\textup{ad}\,g)^2\Lambda(1,\textup{sym}^2\,\phi)}.
  \end{align}
  Similarly,
  \begin{align}\label{eqn:Watson-Ichino E}
    \left|\left<|g|^2,E_\infty\left(\cdot,\frac{1}{2}+i\tau\right)\right>_q\right|^2
    =\frac{1}{4q}
    \left|\frac{\Lambda^q(\frac{1}{2}+i\tau)\Lambda(\frac{1}{2}+i\tau,\textup{ad}\,g)}
    {\Lambda(1,\textup{ad}\,g)\Lambda^q(1+2i\tau)}\right|^2.
  \end{align}
\end{lemma}
Recall we have
\begin{align}\label{eqn:ad g}
  L_\infty(s,\textup{ad}g)=\pi^{-\frac{3s}{2}}\Gamma\left(\frac{s}{2}\right)
  \Gamma\left(\frac{s}{2}+it_g\right)\Gamma\left(\frac{s}{2}-it_g\right),
\end{align}
and
\begin{align}\label{eqn:sym g}
  L_\infty(s,\textup{sym}^2\phi)=\pi^{-\frac{3s}{2}}\Gamma\left(\frac{s}{2}\right)
  \Gamma\left(\frac{s}{2}+it_\phi\right)\Gamma\left(\frac{s}{2}-it_\phi\right).
\end{align}
It is obvious that \eqref{eqn:Watson-Ichino cusp q} vanishes if $\phi$ is odd, so
it is natural to assume $\phi$ is even, in which case
we have
\begin{align}\label{eqn:phi adg}
  L_\infty(s,\phi\times \textup{ad}g)=\pi^{-3s}\prod_{\sigma_1=\pm1}\prod_{\sigma_2=\pm1}
  \Gamma\left(\frac{s+i\sigma_1t_\phi}{2}\right)\Gamma\left(\frac{s+i\sigma_1t_\phi+2i\sigma_2t_g}{2}\right)
\end{align}
and
\begin{align}\label{eqn:phi}
  L_\infty(s,\phi)=\pi^{-s}\Gamma\left(\frac{s+it_\phi}{2}\right)\Gamma\left(\frac{s-it_\phi}{2}\right).
\end{align}
Moreover, since we are working on the spectral aspect, we do not need to care
the arithmetic conductors of the $L$-functions in \eqref{eqn:Watson-Ichino cusp q}
and \eqref{eqn:Watson-Ichino E} which depend on $q$ and $q_1$.

By the Rankin-Selberg unfolding method, we have
(note that $\phi\in\mathcal{B}_0^*(q_1)$)
\begin{align*}
  L(1,\textup{ad}\,g)\asymp \frac{|\rho_g(1)|^2}{\cosh(\pi t_g)}, \quad
  L(1,\textup{sym}^2\,\phi)\asymp \frac{|\rho_\phi(1)|^2}{\cosh(\pi t_{\phi})},
\end{align*}
which give us
\begin{align}\label{eqn:value at 1}
  \frac{1}{L(1,\textup{ad}\,g)}\ll_q (1+|t_g|)^\varepsilon,\quad
  \frac{1}{L(1,\textup{sym}^2\,\phi)}\ll_{q_1} (1+|t_\phi|)^\varepsilon,
\end{align}
by using \eqref{eqn:Fc3}.
We also have the well-known bound (see \cite[3.11.7]{Tit86})
\begin{align}\label{eqn:zeta value at 1}
  \frac{1}{\zeta(1+it)}\ll \log (|t|+3).
\end{align}


\subsection{Summation formulas}\label{subsec:Summation formulas}
We first state the Kuznetsov trace formula.
\begin{lemma}\label{lem:Kuznetsov}
  $\mathrm{(}$\cite[Proposition 2.1]{Mic04}, \cite[Lemma 2.3]{Blomer05}$\mathrm{)}$
  Let $\delta>0$, and let $h$ be a function that is even, holomorphic in the horizontal strip
  $|\Im (t)|\leq \frac{1}{4}+\delta$, and satisfies $h(t)\ll (1+|t|)^{-2-\delta}$.
  Then, for $m$, $n\in\mathbb{N}$,
  \begin{multline*}
    \sqrt{mn}\sum_{\phi\in \mathcal{B}(q,\psi)}\frac{h(t_{\phi})}{\cosh(\pi t_\phi)}
    \overline{\rho_{\phi}(m)}\rho_{\phi}(n)
    +\sqrt{mn}\sum_{\mathfrak{a}\,\mathrm{singular}}
    \frac{1}{4\pi}\int_{\mathbb{R}}
    \frac{h(t)}{\cosh(\pi t)}\overline{\rho_{\mathfrak{a}}(m,t)}\rho_{\mathfrak{a}}(n,t)dt
    \\
    =\frac{\delta_{m,n}}{4\pi^2}\int_{\mathbb{R}}t\tanh(\pi t)h(t)dt
    +\sum_{c\equiv0(q)}\frac{1}{c}S_\psi(m,n;c)H\left(\frac{4\pi\sqrt{mn}}{c}\right),
  \end{multline*}
  where
  \begin{align*}
    H(x)=\frac{i}{2\pi}\int_{\mathbb{R}}\frac{J_{2it}(x)}{\cosh(\pi t)}th(t)dt,
  \end{align*}
  and
  \begin{align*}
    S_\psi(m,n;c)=\underset{d(\mod c)}{\sum\nolimits^*}
    \bar{\psi}(d)e\left(\frac{m\bar{d}+nd}{c}\right)
  \end{align*}
  is the Kloosterman sum.
\end{lemma}

Now we turn to the Voronoi summation formula.
\begin{lemma}\label{lem:VSF2}
  $\mathrm{(}$\cite[Theorem A.4]{KowMicVan02}$\mathrm{)}$
  Let $g\in \mathcal{B}^*(q,\psi)$ with the parity $\epsilon_g$, and
  let $W(x)$ be a a smooth compactly supported function on $(0,\infty)$.
  Then, for $q|c$ and $(d,c)=1$, we have
  \begin{equation}\label{eqn:VSF2}
    \sum_{n\geq1} \lambda_g(n) e\left(\frac{dn}{c}\right) W(n)
    = \frac{\bar{\psi}(d)}{c} \sum_{\pm}
    \sum_{n\geq1} \lambda_g(n) e\left(\mp\frac{\bar{d} n}{c}\right) W^\pm \left(\frac{n}{c^2}\right)
  \end{equation}
  where
  \begin{equation}\label{eqn:W+}
    W^+ (y) =   \frac{\pi i}{\sinh(\pi t_g)} \int_{0}^{\infty} W(x) (J_{2it_g}(4\pi\sqrt{yx})-J_{-2it_g}(4\pi\sqrt{yx})) d x,
  \end{equation}
  and
  \begin{equation}\label{eqn:W-}
    W^- (y) =  4 \epsilon_g \cosh(\pi t_g) \int_{0}^{\infty} W(x) K_{2it_g}(4\pi\sqrt{yx})
    dx.
  \end{equation}
\end{lemma}

\section{Analytic preliminaries}\label{sec:analytic preliminaries}

\subsection{Smooth weight functions}\label{subsec:smooth weight functions}

We use the notation as in \cite[\S 2.4]{BloJanNel}.
Let $T$ be the large parameter in Theorem \ref{thm:moment result}.
We introduce the abbreviation
\begin{align*}
  A_1\preccurlyeq A_2\Longleftrightarrow A_1\ll_\varepsilon T^\varepsilon A_2.
\end{align*}
A smooth function $V:\mathbb{R}^n\rightarrow \mathbb{C}$
is called flat if
\[
  x_1^{j_1} \cdots x_n^{j_n} V^{(j_1,\ldots,j_n)} (x_1,\ldots,x_n)
   \preccurlyeq_{\mathbf{j}} 1
\]
for all $\mathbf{j}\in\mathbb{Z}^n_{\geq0}$.
It is obvious that if $V$ is flat, then $\exp(iV)$ is also flat.
We call $V$ nice if it is flat and has compact support in $(0,\infty)^n$.
From now on, we denote by $V$ a nice function in one or more variables,
and we will redefine them from line to line to suit our needs.

We may separate variables in $V(x_1,\ldots,x_n)$ by the Mellin transform.
Precisely, let $V(x_1)\cdots V(x_n)$ be a nice function which is $1$
on the support of $V$. Then by the Mellin inversion, one has
\begin{multline*}
  V(x_1,\ldots,x_n)=V(x_1,\ldots,x_n)V(x_1)\cdots V(x_n)
  \\
  =\int_{\Re s_1=0}\cdots \int_{\Re s_n=0}\hat{V}(s_1,\ldots,s_n)
  \left(V(x_1)\cdots V(x_n)x_1^{-s_1}\cdots x_n^{-s_n}\right)
  \frac{ds_1\cdots ds_n}{(2\pi i)^n}.
\end{multline*}
We can truncate the vertical integrals at $|\Im s_i|\preccurlyeq 1$
at the cost of a negligible error.
We will often use this technique to separate variables without explicit mention.

\subsection{Oscillatory integrals}

We will frequently use the stationary phase method.
So we quote the following lemmas.
\begin{lemma}\label{lemma:stationary_phase BKY}
Let $Y\geq 1$, $X$, $P$, $U$, $S>0$, and suppose that $w$ is a smooth function with support
on $[\alpha,\beta]$, satisfying
\begin{align*}
  w^{(j)}(t)\ll_j XU^{-j}.
\end{align*}
Suppose $h$ is a smooth function on $[\alpha,\beta]$ such that
\begin{align*}
 |h'(t)|\gg S
 \end{align*}
  for some $S>0$, and
 \begin{align*}
    h^{(j)}\ll_j YP^{-j},\quad for \ j=2,3,\cdots.
  \end{align*}
 Then
 \begin{align*}
 \int_{\mathbb{R}}w(t)e^{ih(t)}dt\ll_A (\beta-\alpha)X[(PS/\sqrt{Y})^{-A}+(SU)^{-A}].
 \end{align*}
\end{lemma}
Recently, this result has been extended by (see \cite{KirPetYou19}).
The following lemma is due to in \cite[\S 3, Main Theorem]{KirPetYou19}.
For our convenience, we use the statement as in \cite[Lemma 4]{BloJanNel}.
\begin{lemma}\label{lem:main theorem in KPY}
  Let $T$ be a large parameter. Let $V(t_1,\ldots, t_d)$ be a flat function in the sense of
  \S \ref{subsec:smooth weight functions} with support in $\times_{j=1}^d[c_{1j},c_{2j}]$
  for some fixed intervals $[c_{1j},c_{2j}]\subseteq \mathbb{R}$ not containing $0$.
  Let $X_1$, $X_2$, $\ldots$, $X_d>0$, $Y\geq T^\varepsilon$.
  Write $\mathscr{S}=\times_{j=1}^d[c_{1j}X_j,c_{2j}X_j]\subseteq \mathbb{R}^d$.
  Let $\phi:\mathbb{R}^d\rightarrow \mathbb{R}$ be a smooth function satisfying the derivative
  upper bounds
  \begin{align*}
    \phi^{(j_1,j_2,\ldots,j_d)}(t_1,t_2,\ldots,t_d)\preccurlyeq Y\prod_{i=1}^dX_i^{-j_i}
  \end{align*}
  for $\mathbf{j}\in \mathbb{N}_0^d$ and $(t_1,\ldots,t_d)\in\mathscr{S}$,
  as well as the following second derivative lower bound in the first variable:
  \begin{align*}
    \phi^{(2,0,\ldots,0)}(t_1,t_2,\ldots,t_d)\gg YX_1^{-2}.
  \end{align*}
  Suppose that there exists $t^*=t^*(t_2,\ldots,t_d)$ such that
  $\phi^{(1,0,\ldots,0)}(t^*,t_2,\ldots,t_d)=0$.
  Then
   \begin{align*}
     \int_{\mathbb{R}}V\left(\frac{t_1}{X_1},\ldots,\frac{t_d}{X_d}\right)
     e^{i\phi(t_1,\ldots,t_d)}dt_1=\frac{X_1}{Y^{1/2}}e^{i\phi(t^*,t_2,\ldots,t_d)}
     W\left(\frac{t_2}{X_2},\ldots,\frac{t_d}{X_d}\right)+O(T^{-B})
   \end{align*}
   for a flat function $W=W_B$ for every $B>0$ with support in $\times_{j=2}^d[c_{1j},c_{2j}]$.
\end{lemma}


We will need to estimate integrals of the form
\begin{align}\label{eqn:typical integral in BJN}
  \int_{\mathbb{R}}V\left(\frac{x}{M}\right)
  e(\alpha x^{\frac{1}{2}}+\beta x^{-\frac{1}{2}}) dx
\end{align}
for certain $\alpha$, $\beta\in\mathbb{R}$ satisfying
\begin{align}\label{eqn:condition of sub5.4 in BJN}
  |\alpha|M^{\frac{1}{2}} + |\beta|M^{-\frac{1}{2}} \gg M^\varepsilon.
\end{align}
This has been explicitly computed in \cite[\S 5.4]{BloJanNel} by using
Lemma \ref{lemma:stationary_phase BKY} and Lemma \ref{lem:main theorem in KPY}.
For our convenience, we state it as the following lemma.
\begin{lemma}\label{lem:result of sub5.4 in BJN}
  Under the condition of \eqref{eqn:condition of sub5.4 in BJN}, the integral
  \eqref{eqn:typical integral in BJN} is negligible unless $\frac{\beta}{\alpha}\asymp M$,
  in which case we introduce two parameters $K_1$ and $K_2$, such that
   one can restrict to dyadic ranges $\alpha\asymp K_1$, $\beta\asymp K_2$
  and also possibly restrict the support of $V$ to a neighbourhood of $t^*=\frac{\beta}{\alpha}$.
  Then, we have
  \begin{align*}
    \eqref{eqn:typical integral in BJN}
    =\frac{M^{\frac{5}{4}}}{|\beta|^{\frac{1}{2}}}V\left(\frac{\alpha}{K_1}\right)
    V\left(\frac{\beta}{K_2}\right)V\left(\frac{\beta/\alpha}{M}\right)
    e(2\sgn(\alpha)\sqrt{\alpha\beta})+O(M^{-1000}),
  \end{align*}
  with different functions $V$.
\end{lemma}

\subsection{Bessel functions}

In this subsection, we collect some results for Bessel functions.
By \cite[7.13.2(17)]{EMOT}, we have
\begin{align}\label{eqn:J2it}
\begin{split}
  &\frac{J_{2it}(2x)}{\cosh(\pi t)}=\sum_{\pm}e^{\pm 2i\omega(x,t)}
  \frac{f_{A}^\pm(x,t)}{x^{1/2}+|t|^{1/2}}
  +O_A((x+|t|)^{-A}),
  \\
  &\omega(x,t)=|t|\textup{arcsinh}\frac{|t|}{x}-\sqrt{t^2+x^2},
\end{split}
\end{align}
for $t\in\mathbb{R}$, $|t|>1$ and $x>0$,
where for any fixed $A>0$ the function $f_A^\pm$ is flat.
It is easy to compute that
\begin{align}\label{eqn:properties of omega}
  \begin{split}
     \frac{\partial}{\partial t}\omega(x,t)=\textup{arcsinh}\frac{t}{x},
     \quad
     \frac{\partial^2}{\partial t^2}\omega(x,t)=\frac{1}{\sqrt{x^2+t^2}},
     \quad
     \frac{\partial}{\partial x}\omega(x,t)=-\frac{\sqrt{x^2+t^2}}{x},
  \end{split}
\end{align}
and for $t\ll x^{\frac{3}{4}-\varepsilon}$, one has
\begin{align}\label{eqn:Taylor of omega}
  \omega(x,t)=-x+\frac{t^2}{2x}+O\left(\frac{t^4}{x^3}\right).
\end{align}

We proceed with some results of the $K$-Bessel function.
By \cite{Bal}, we have
\begin{align}\label{eqn:bound K Bessel}
  \cosh(\pi t)K_{2it}(2x)\ll
  \begin{cases}
     t^{-\frac{1}{4}}(t-x)^{-\frac{1}{4}}, \quad &0<x<t-C_0t^{\frac{1}{3}}, \\
     t^{-\frac{1}{3}}, \quad &|x-t|\leq C_0t^{\frac{1}{3}}, \\
     x^{-\frac{1}{4}}(x-t)^{-\frac{1}{4}}
     \exp(-c_0(\frac{x}{t})^{\frac{3}{2}}(\frac{x-t}{t^{1/3}})^{\frac{3}{2}}), \quad &x>t+C_0t^{\frac{1}{3}},
  \end{cases}
\end{align}
for $t\geq 1$. Here, $c_0$ and $C_0$ are two positive constants.
Furthermore, one has the integration representation \cite[8.432-4]{GraRyz}
\begin{align}\label{eqn:integration rep of K}
  2\cosh(\pi t)K_{2it}(2x)
  =\int_{-\infty}^\infty \cos(2x\sinh v)\exp(2itv)dv,
\end{align}
for $t\in\mathbb{R}$ and $x>0$.


\section{Reduction of Theorem \ref{thm:moment result}}\label{sec:reduction}
We only focus on the proof of \eqref{eqn:main moment 1}.
Actually, we will see that \eqref{eqn:main moment 2}
is the special case of \eqref{eqn:main moment 1}.

Let
\begin{align*}
  h(t)=e^{-\frac{(t-T)^2}{H^2}}+e^{-\frac{(t+T)^2}{H^2}},
\end{align*}
and recall
\begin{align}\label{eqn:parameters of tf pm tg}
  T_0=T+t_g,\quad L=|T-t_g|, \quad H=o(L).
\end{align}
To prove \eqref{eqn:main moment 1}, by using Lemma \ref{lemma:AFE-LiYoung}
together with by \eqref{eqn:Fc1}, \eqref{eqn:Fc3},
it is enough to get
\begin{align}\label{eqn:1st reduce}
  \sum_{f\in \mathcal{B}^*(D,\chi)}\frac{h(t_f)}{\cosh (\pi t_f)}
  \left|\sum_{n\geq 1}\sqrt{n}\rho_f(n)\lambda_g(n)V\left(\frac{n}{N}\right)\right|^2
  \preccurlyeq NT_0L^{\frac{1}{2}}||g||_4^2+NTH,
\end{align}
for all $N\preccurlyeq T_0L$.
Here we have made dyadic decomposition and absorbed the factor $n^{-iw}$ to the
smooth weight function $V$.
Then, since the LHS of \eqref{eqn:1st reduce} is positive,
we may enlarge $\mathcal{B}^*(q,\chi)$ to $\mathcal{B}([q,D],\chi)$
and use the Kuznetsov trace formula (Lemma \ref{lem:Kuznetsov}).
Note that the part of the continuous spectrum is positive and the diagnoal term is acceptable
(bounded by $\preccurlyeq NTH$),
it is sufficient to consider
\begin{align*}
  \sum_{n\geq 1}\sum_{m\geq 1}\lambda_g(n)\lambda_g(m)V\left(\frac{n}{N}\right)
  V\left(\frac{m}{N}\right)
  \sum_{c\equiv 0(\mod([q,D]))}\frac{S_{\chi}(m,n;c)}{c}\int_{\mathbb{R}}e^{-\frac{(t-T)^2}{H^2}}
  \frac{tJ_{2it}(2x)}{\cosh(\pi t)}\textup{d}t,
\end{align*}
where $x=\frac{2\pi\sqrt{nm}}{c}$.
Now we use the treatment in \cite[p.39]{BloJanNel}
to truncate $c$.
We can first truncate $c$ by some large power of $T$ by shifting the contour
of the $t$-integral (see \cite[p.75]{JutMot05}).
So the integral can be smoothly truncated at $t\in[T-HT^\varepsilon, T+HT^\varepsilon]$.
Then, by using by applying \eqref{eqn:J2it},
\eqref{eqn:properties of omega} and Lemma \ref{lemma:stationary_phase BKY}
with $S=\min(1,\frac{T}{x})$, $P=Y=T+x$, $U=H$,
one sees that the integral is negligible unless $x\gg T^{1-\varepsilon}H$,
in which case we have
$c\preccurlyeq \frac{N}{TH}$
and
$t\ll T_0\ll x^{\frac{3}{4}-\varepsilon}$.
By using \eqref{eqn:Taylor of omega}, we obtain
\begin{align*}
  \frac{J_{2it}(2x)}{\cosh(\pi t)}=x^{-1/2}\sum_{\pm}e^{\pm2i(-x+\frac{t^2}{2x})}F^\pm(x,t)+O(T^{-1000})
\end{align*}
for a flat function $F^\pm$.


Hence, by smoothing the parameter $c$, it suffices to get
\begin{multline}\label{re goal}
  \frac{TH}{(CN)^{1/2}}\sum_{n\geq 1}\sum_{m\geq 1}\lambda_g(n)\lambda_g(m)
  V\left(\frac{n}{N}\right)V\left(\frac{m}{N}\right)
  \sum_{c\equiv 0(\mod([q,D]))}S_{\chi}(m,n;c)V\left(\frac{c}{C}\right)
  \\
  \cdot\sum_{\sigma_1=\pm 1}
  e\left(\sigma_1\frac{2\sqrt{mn}}{c}-\sigma_1\frac{t^2c}{4\pi^2\sqrt{mn}}\right)
  \preccurlyeq NT_0L^{\frac{1}{2}}||g||_4^2+NTH,
\end{multline}
uniformly in $t\in [T-HT^\varepsilon,T+HT^\varepsilon]$ and $C\preccurlyeq \frac{N}{TH}$.
We open the Kloosterman sum and use the Voronoi summation formula (Lemma \ref{lem:VSF2})
to the $n$-sum, getting
\begin{multline}\label{Voronoi for n}
  \sum_{n\geq 1}\lambda_g(n)V\left(\frac{n}{N}\right)e\left(\frac{dn}{c}\right)
  e\left(\sigma_1\frac{2\sqrt{mn}}{c}-\sigma_1\frac{t^2c}{4\pi^2\sqrt{mn}}\right)
  \\
  =\frac{\psi(\bar{d})}{c}\sum_\pm\sum_{n\geq 1}
  \lambda_g(n)e\left(\mp\frac{\bar{d}n}{c}\right)
  \int_0^\infty V\left(\frac{x}{N}\right)
  e\left(\sigma_1\frac{2\sqrt{mx}}{c}-\sigma_1\frac{t^2c}{4\pi^2\sqrt{mx}}\right)
  \mathcal{J}^\pm(2y)dx,
\end{multline}
where $y=\frac{2\pi\sqrt{nx}}{c}$ and $\mathcal{J}^\pm$ are defined as
in \eqref{eqn:W+} and \eqref{eqn:W-} with $W(x)$ replaced by
$V(\frac{x}{N})e(\sigma_1\frac{2\sqrt{mx}}{c}-\sigma_1\frac{t^2c}{4\pi^2\sqrt{mx}})$.

\section{The terms related to $\mathcal{J}^+$}
In this section, we deal with the contribution from the $\mathcal{J}^+$-term
in \eqref{Voronoi for n}.
Note that
\begin{align*}
  \mathcal{J}^+(2y)=\pi i\frac{\cosh (\pi t_g)}{\sinh(\pi t_g)}\frac{J_{2it_g}(2y)-J_{-2it_g}(2y)}{\cosh(\pi t_g)}.
\end{align*}
We will only deal with the $J_{2it_g}$-term,
since the $J_{-2it_g}$-term is similar ($J_{-2it_g}(2y)=\overline{J_{2it_g}(2y)}$).
By \eqref{eqn:J2it}, the $x$-integral in \eqref{Voronoi for n} is led to
\begin{align}\label{eqn:reduce x-integral}
  \int_0^\infty V\left(\frac{x}{N}\right)
  e\left(\sigma_1\frac{2\sqrt{mx}}{c}-\sigma_1\frac{t^2c}{4\pi^2\sqrt{mx}}
  \pm\frac{1}{\pi}\omega(y,t_g)\right)\frac{f^{\pm}_A(y,t_g)}{y^{1/2}+t_g^{1/2}}dx.
\end{align}
By \eqref{eqn:properties of omega}, we have
\begin{align*}
  x\frac{\partial \omega(y,t_g)}{\partial x}\asymp \frac{\sqrt{nN}}{c}+t_g,
\end{align*}
and
\begin{align*}
  x^j\frac{\partial^j \omega(y,t_g)}{\partial x^j}\preccurlyeq \frac{\sqrt{nN}}{c}+t_g.
\end{align*}
Note that $f^{\pm}_A(y,t_g)$ is also a flat function for $x$
and $t_g\ll (TH)^{3/4-\varepsilon}=o(\frac{N}{C})$.
So, we apply \eqref{eqn:properties of omega} and
Lemma \ref{lemma:stationary_phase BKY} with

\begin{align*}
  X=\left(\frac{\sqrt{nN}}{C}+t_g\right)^{-\frac{1}{2}}T^\varepsilon,\quad P=U=N,
  \quad Y=\frac{N}{C}+\frac{\sqrt{nN}}{C},
  \quad S=\frac{Y}{N}
\end{align*}
seeing that the integral is negligible unless $n\asymp N$, in which case,
we use \eqref{eqn:Taylor of omega} and get
\begin{align}\label{eqn:app J2it}
  \frac{J_{2it_g}(2y)}{\cosh(\pi t_g)}=
  y^{-1/2}\sum_{\pm}e^{\pm2i(-y+\frac{t_g^2}{2u})}F^\pm(y,t_g)+O(T^{-1000}),
\end{align}
where $F^\pm$ is flat.
Thus \eqref{eqn:reduce x-integral} is reduced to
\begin{align}\label{x integral}
  \frac{C^{\frac{1}{2}}}{N^{\frac{1}{2}}}\int_{\mathbb{R}}V\left(\frac{x}{N}\right)
  e(\beta_1x^{-1/2})e(\alpha x^{1/2}+\beta x^{-1/2})\textup{d}x,
\end{align}
where
\begin{align*}
  \alpha=\frac{2\sigma_1\sqrt{m}+2\sigma_2\sqrt{n}}{c}, \quad
  \beta=-\frac{\sigma_1c}{4\pi^2\sqrt{m}}(t^2-t_g^2), \quad
  \beta_1=-\frac{t_g^2c}{4\pi^2}\left(\frac{\sigma_1}{\sqrt{m}}+\frac{\sigma_2}{\sqrt{n}}\right),
\end{align*}
with $\sigma_1$, $\sigma_2=\pm1$.

If $|\alpha| N^{1/2}+|\beta|N^{-1/2}\preccurlyeq 1$,
then the contribution is admissible.
Actually, we can first get
$c\asymp C\preccurlyeq\frac{N}{T_0L}\preccurlyeq 1$
from $|\beta|N^{-1/2}\preccurlyeq 1$.
Then, by using this and $|\alpha| N^{1/2}\preccurlyeq 1$,
we get $\sigma_1=-\sigma_2$
and $|m-n|\preccurlyeq C\preccurlyeq 1$.
Thus the contribution to the LHS of \eqref{re goal} is
\begin{align*}
  &\preccurlyeq \frac{TH}{(CN)^{1/2}}\sum_{m\asymp N}|\lambda_g(m)|
  \sum_{\substack{n\asymp N\\ |m-n|\preccurlyeq 1}}|\lambda_g(n)|
  \frac{N^{1/2}}{C^{1/2}}\sum_{c\asymp C}|S_{\chi\psi}((m-n),0;c)|
  \\
  &\preccurlyeq TH\sum_{\substack{m,n\asymp N\\ |m-n|\preccurlyeq 1}}
  (\lambda_g^2(m)+\lambda_g^2(n))
  \\
  &\preccurlyeq NTH,
\end{align*}
which is acceptable.

From now on, we assume that $|\alpha| N^{1/2}+|\beta|N^{-1/2}\gg T^\varepsilon$.
It is clear that
\begin{align*}
  \beta_1x^{-\frac{1}{2}}(\alpha x^{\frac{1}{2}})^{-1}
  \ll \frac{T_0^2C^2}{N^2}\preccurlyeq (TH)^{-\frac{1}{2}},
\end{align*}
and
\begin{align*}
  \beta x^{-\frac{1}{2}}\left(\frac{\sqrt{mx}}{c}\right)^{-1}
  \ll \frac{T_0^2C^2}{N^2}\preccurlyeq (TH)^{-\frac{1}{2}}.
\end{align*}
By using Lemma \ref{lemma:stationary_phase BKY}
with $S=|\alpha|^{1/2}N^{-1/2}+|\beta|N^{-3/2}$,
$P=U=N$, $X=T^\varepsilon$ and $Y=NS$, the $x$-integral in \eqref{x integral}
is negligible unless
$\sigma_1=-\sigma_2$ and
$|\alpha x^{1/2}|\asymp |(\alpha x^{1/2}+\beta_1x^{-1/2})|\asymp |\beta x^{-1/2}|$,
in which case we have
\begin{align*}
  r:=|n-m|\asymp R:=\frac{T_0LC^2}{N}\preccurlyeq N(TH)^{-\frac{1}{2}-\varepsilon}.
\end{align*}
Actually we can also get $\sigma_3:=\sgn(n-m)=\sgn(t-t_g)$.
By inserting a nice function $V(\frac{r}{R})$ to smooth the parameter $r$,
we see that
\begin{align*}
  e(\beta_1x^{-1/2})=e\left(-\frac{\sigma_1\sigma_3}{4\pi^2}
  \frac{t_g^2cr}{\sqrt{m(m+\sigma_3r)}(\sqrt{m}+\sqrt{m+\sigma_3r})}x^{-\frac{1}{2}}\right)
\end{align*}
can be absorbed into the flat functions
$V(\frac{m}{N})V(\frac{r}{R})V(\frac{c}{C})V(\frac{x}{N})$, since
\begin{align}\label{eqn:bound beta1}
  \beta_1\ll \frac{t_g^2CR}{N^{\frac{3}{2}}}\ll \frac{t_g^2T_0LC^3}{N^{\frac{5}{2}}}
  \preccurlyeq \frac{N^{\frac{1}{2}}T_0^4}{(TH)^3}\preccurlyeq N^{\frac{1}{2}}.
\end{align}
Consequently, we may replace \eqref{eqn:reduce x-integral} (or \eqref{x integral}) by
\begin{align*}
  \frac{C^{\frac{1}{2}}}{N^{\frac{1}{2}}}
  \int_{\mathbb{R}}V\left(\frac{x}{N}\right)
  e(\alpha x^{1/2}+\beta x^{-1/2})\textup{d}x.
\end{align*}
By Lemma \ref{lem:result of sub5.4 in BJN},
the contribution of the $\mathcal{J}^+$-term to \eqref{Voronoi for n} is
\begin{align}\label{eqn:Voronoi J contribution}
  \frac{\bar{\psi}(d)N}{(T_0L)^{1/2}C}
  \sum_{n}V\left(\frac{r}{R}\right)\lambda_g(n)e\left(-\frac{\bar{d}n}{c}\right)
  e\left(-2\sigma\left(\frac{|t^2-t_g^2|}{2\pi^2 \sqrt{m}}\right)^{1/2}
  |\sqrt{n}-\sqrt{m}|^{1/2}\right),
\end{align}
where $\sigma=\sigma_1\sigma_3$.
Therefore, to estimate the contribution from $\mathcal{J}^+$ to the LHS
of \eqref{re goal},
it suffices to get the following bound
\begin{multline}\label{eqn:rere goal}
  \frac{TH}{C^{3/2}}\left(\frac{N}{T_0L}\right)^{1/2}
  \sum_{c\equiv0(\mod[q,D])}V\left(\frac{c}{C}\right)
  \sum_{r\geq 1}
  G_{\chi\psi}(r,c)V\left(\frac{r}{R}\right)
  \\
  \cdot
  \underset{\substack{m,n\geq 1\\{|n-m|=r}}}{\sum\sum}
  \lambda_g(n)\lambda_g(m)
  V\left(\frac{n}{N}\right)V\left(\frac{m}{N}\right)
  e\left(-2\sigma\left(\frac{|t^2-t_g^2|}{2\pi^2 \sqrt{m}}\right)^{1/2}
  |\sqrt{n}-\sqrt{m}|^{1/2}\right)
  \\
  \preccurlyeq NT_0L^{\frac{1}{2}}||g||^2_4+NTH,
\end{multline}
where (note that $\chi(-1)=\psi(-1)=1$)
\begin{align*}
  G_{\chi\psi}(r,c)=
  \underset{d(c)}{\sum\nolimits^*}\overline{\chi\psi}(d)
  e\left(-\sigma_3\frac{\bar{d}r}{c}\right)
  =\underset{d(c)}{\sum\nolimits^*}\chi\psi(d)e\left(\frac{dr}{c}\right).
\end{align*}
For our convenience, we denote by $L_0:=|t-t_g|$.
So, by $t\in [T-HT^\varepsilon,T+HT^\varepsilon]$
and $H=o(L)$, one has $L_0\asymp L$.

If $t>t_g$ (which means $n>m$),
by
\begin{align*}
  (\sqrt{n}-\sqrt{m})^{1/2}=m^{1/4}\left(\frac{r}{2m}\right)^{1/2}
  \left(1-\frac{r}{8m}+\cdots\right),
\end{align*}
we may replace $\sqrt{\sqrt{n}-\sqrt{m}}$ in \eqref{eqn:Voronoi J contribution}
by $m^{1/4}(\frac{r}{2m})^{1/2}$,
with the error term absorbed into the flat functions
$V(\frac{m}{N})V(\frac{r}{R})$.
In fact, we have
\begin{align*}
  m^{1/4}\left(\frac{r}{m}\right)^{3/2}\left(\frac{|t^2-t_g^2|}{m^{1/2}}\right)^{1/2}
  \ll \frac{R^{3/2}}{N^{3/2}}(T_0L)^{1/2}\preccurlyeq
   T^{-\varepsilon}.
\end{align*}
Then, just for our convenience (not necessary),
we replace $|t^2-t_g^2|$ in \eqref{eqn:rere goal}
by $T_0L_0$ up to a flat function,
since
\begin{align*}
  |t^2-t_g^2|^\frac{1}{2} =(T_0L_0)^{\frac{1}{2}}
  \left(1+\frac{t-T}{2T_0}+\cdots\right),
\end{align*}
and
\begin{align*}
  \frac{(T_0L_0)^{1/2}}{m^{1/4}}\frac{|t-T|}{T_0}m^{1/4}\left(\frac{r}{m}\right)^{1/2}
  \preccurlyeq \left(\frac{LR}{T_0N}\right)^{1/2}H\preccurlyeq \frac{LCH}{N}
  \preccurlyeq \frac{L}{T}.
\end{align*}
This allows us to rewrite the second line of \eqref{eqn:rere goal} as
\begin{align*}
  \sum_{m\geq 1}\lambda_g(m)\lambda_g(m+r)
  V\left(\frac{m}{N}\right)
  e\left(-\sigma\frac{(T_0L_0r)^{\frac{1}{2}}}{\pi m^{\frac{1}{2}}}\right).
\end{align*}
For $t<t_g$, we have $m>n$. Similarly, the second line of \eqref{eqn:rere goal}
can be recast as
\begin{align}\label{eqn:t<tg}
  \sum_{n\geq 1}\lambda_g(n)\lambda_g(n+r)
  V\left(\frac{n}{N}\right)
  e\left(-\sigma\frac{(T_0L_0r)^{\frac{1}{2}}}{\pi (n(n+r))^{\frac{1}{4}}}\right).
\end{align}
By
$(n+r)^{-1/4}=n^{-1/4}(1-\frac{r}{4n}+\frac{5r^2}{32n^{2}}+\cdots)$,
we have
\begin{align*}
  \frac{(T_0L_0r)^{\frac{1}{2}}}{(n(n+r))^{\frac{1}{4}}}
  -\frac{(T_0L_0r)^{\frac{1}{2}}}{n^{\frac{1}{2}}}
  \ll \frac{(T_0L_0)^{\frac{1}{2}}R^{\frac{3}{2}}}{N^{\frac{3}{2}}}
  \preccurlyeq T^{-\varepsilon}.
\end{align*}
So \eqref{eqn:t<tg} can be replaced by
\begin{align*}
   \sum_{n\geq 1}\sum_{r\geq 1}\lambda_g(n)\lambda_g(n+r)
  V\left(\frac{n}{N}\right)
  e\left(-\sigma\frac{(T_0L_0r)^{\frac{1}{2}}}{\pi n^{\frac{1}{2}}}\right),
\end{align*}
with the error term absorbed into the flat functions
$V(\frac{n}{N})V(\frac{r}{R})$.

Therefore, to get \eqref{eqn:rere goal},
we are reduced to showing
\begin{multline}\label{rere goal Maass}
  \frac{TH}{C^{3/2}}\left(\frac{N}{T_0L}\right)^{1/2}
  \sum_{r\geq 1}P(r,N)V\left(\frac{r}{R}\right)
  \sum_{c\equiv0(\mod[q,D])}V\left(\frac{c}{C}\right)
  G_{\chi\psi}(r,c)
  \\
  \preccurlyeq NT_0L^{\frac{1}{2}}||g||^2_4+NTH,
\end{multline}
where
\begin{align}\label{eqn:SVS}
  P(r,N)=\sum_{m\geq 1}\lambda_g(m)\lambda_g(m+r)V\left(\frac{m}{N}\right)
  e\left(-\sigma\frac{(T_0L_0r)^{1/2}}{\pi m^{1/2}}\right).
\end{align}
Note that we have reindexed $n$ and $m$ when $t<t_g$.

If $R\preccurlyeq 1$, then it is clear that
$C^2\preccurlyeq \frac{N}{T_0L}\preccurlyeq 1$.
Thus the contribution of $R\preccurlyeq 1$ in \eqref{rere goal Maass}
is
\begin{align*}
  \preccurlyeq NTH\left(\frac{N}{T_0L}\right)^{1/2}
  \preccurlyeq NTH,
\end{align*}
which is acceptable.

From now on, we  assume $R\gg T^\varepsilon$.
Let $c=[q,D]^kc_k$, where $(c_k,qD)=1$.
Then, it is obvious that $k\ll \log T$.
Note that in the case of $q=D=1$, we take $k=1$.
By applying the change of variable $d=d_1[q,D]^k+d_2c_k$ with
$(d_1,c_k)=(d_2,qD)=1$, one gets
\begin{align*}
  G_{\chi\psi}(r,c)&=
  \underset{d_2([q,D]^k)}{\sum\nolimits^*}\chi\psi(d_2c_k)e\left(\frac{d_2r}{[q,D]^k}\right)
  \underset{d_1(c_k)}{\sum\nolimits^*}e\left(\frac{d_1r}{c_k}\right)
  \\
  &=\chi\psi(c_k)G_{\chi\psi}(r,[q,D]^k)\sum_{d|(c_k,r)}\mu\left(\frac{c_k}{d}\right)d.
\end{align*}
So we can replace the $c$-sum
$\sum_{c\equiv0(\mod[q,D])}V\left(\frac{c}{C}\right)G_{\chi\psi}(r,c)$
in \eqref{rere goal Maass} by
\begin{align}\label{eqn:replace the c sum}
  \sum_{1\leq k\ll \log T}\sum_{\substack{c_k\geq 1\\(c_k,qD)=1}}V\left(\frac{c_k}{C_k}\right)
  \chi\psi(c_k)G_{\chi\psi}(r,[q,D]^k)\sum_{d|(c_k,r)}\mu\left(\frac{c_k}{d}\right)d,
\end{align}
where $[q,D]^kC_k=C$.

Now we are ready to treat $P(r,N)$. Let $N_0=NT^{-\varepsilon}$.
Then, by applying $\lambda_g(m)\lambda_g(m+r)=\lambda_g(-m)\lambda_g(-m-r)$,
and changing variables $-m\rightarrow m\rightarrow m+r$, we get
\begin{multline}\label{eqn:nead to be real}
  \sum_{m<-N_0}\lambda_g(m)\lambda_g(m+r)
  V\left(\frac{|2m+r|-r}{2N}\right)
  e\left(-\sigma\frac{(T_0L_0r)^{\frac{1}{2}}}{\pi (\frac{|2m+r|-r}{2})^{\frac{1}{2}}}\right)
  \\
  =\sum_{m>N_0}\lambda_g(m+r)\lambda_g(m)
  V\left(\frac{m}{N}\right)
  e\left(-\sigma\frac{(T_0L_0r)^{\frac{1}{2}}}{\pi m^{\frac{1}{2}}}\right),
\end{multline}
which implies that
\begin{align*}
  2P(r,N)
  =\sum_{|m|>N_0}\lambda_g(m)\lambda_g(m+r)V\left(\frac{|2m+r|-r}{2N}\right)
  e\left(-\sigma\frac{(T_0L_0r)^{\frac{1}{2}}}{\pi (\frac{|2m+r|-r}{2})^{\frac{1}{2}}}\right).
\end{align*}
In the above deduction, we need the fact that $\psi$ is real.
In fact, by noting that $P(r,N)=\sum_{m\geq 1}\bar{\lambda}_g(m)\lambda_g(m+r)\{\cdots\}$,
we can only get
\begin{align*}
  P(r,N)
  +\sum_{m\geq 1}\lambda_g(m)\bar{\lambda}_g(m+r)\{\cdots\}
  =\sum_{|m|>N_0}\bar{\lambda}_g(m)\lambda_g(m+r)\{\cdots\},
\end{align*}
from the above evaluation.
Denote by
\begin{align*}
G_r(w)=\left(\frac{2z+\frac{r}{N}}{\sqrt{z(z+\frac{r}{N})}}\right)^{2it_g}
V(z)e\left(-\sigma\frac{(T_0L_0r)^{\frac{1}{2}}}{\pi(Nz)^{\frac{1}{2}}}\right),
\end{align*}
where $w=z+\frac{r}{2N}$.
Let $\tilde{G}_r(s)$ be the Mellin transform of $G_r(w)$.
Then, one has
\begin{align*}
  \tilde{G}_r(s)=\int_{0}^\infty
  \left(\frac{2z+\frac{r}{N}}{\sqrt{z(z+\frac{r}{N})}}\right)^{2it_g}
  V(z)e\left(-\sigma\frac{(T_0L_0r)^{\frac{1}{2}}}
  {\pi (Nz)^{\frac{1}{2}}}\right)\left(z+\frac{r}{2N}\right)^{s-1}\textup{d}z.
\end{align*}
Now let $w=\frac{|2m+r|}{2N}$ and $|m|>N_0$.
Then, we have $w=\frac{|m|+|m+r|}{2N}$ and
$\frac{2z+r/N}{\sqrt{z(z+r/N)}}=\frac{|m|+|m+r|}{\sqrt{|m||m+r|}}.$
Hence, we get
\begin{align*}
  &\sum_{|m|>N_0}\lambda_g(m)\lambda_g(m+r)V\left(\frac{|2m+r|-r}{2N}\right)
  e\left(-\sigma\frac{(T_0L_0r)^{\frac{1}{2}}}{\pi (\frac{|2m+r|-r}{2})^{\frac{1}{2}}}\right)
  \\
  =&\frac{1}{2\pi i}\int_{(2)}\tilde{G}_r(s)\sum_{m>N_0}\lambda_g(m)\lambda_g(m+r)
  \left(\frac{\sqrt{|m||m+r|}}{|m|+|m+r|}\right)^{2it_g}
  \left(\frac{|m|+|m+r|}{2N}\right)^{-s}ds,
\end{align*}
which implies that
\begin{align}\label{eqn:express P by D}
  P(r,N)=\frac{1}{4\pi i}\int_{(2)}(2N)^s\tilde{G}_r(s)D_{g}^{\dag}(s,1,1,r)ds,
\end{align}
where $D_{g}^{\dag}(s,\ell_1,\ell_2,r)$ is the same as $D_{g}(s,\ell_1,\ell_2,r)$
but restricted to $|m|>N_0$.
We can also replace $D_{g}^{\dag}$ by $D_{g}$ in \eqref{eqn:express P by D},
since, if $|m|\leq N_0$, we have
\begin{align*}
  \frac{1}{2\pi i}\int_{(2)}(2N)^s\tilde{G}_r(s)(|m|+|m+r|)^{-s}ds
  =G_r\left(\frac{|m|+|m+r|}{2N}\right)=0.
\end{align*}
However, we would like to use the expression \eqref{eqn:express P by D} here.

\subsection{Treating $\tilde{G}_{r}(s)$}
Let $s=u+iv$ be in a sufficiently large but fixed strip
$-A\leq u\leq A$ (where we may take $A$ to be large).
Then, we have
\begin{align*}
  \tilde{G}_{r}(s)=4^{it_g}\int_0^\infty V(z)\left(z+\frac{r}{2N}\right)^{u-1}
  e(\phi(z)+h(z))\textup{d}z,
\end{align*}
where
\begin{align*}
  \phi(z)=\phi(z,r,v)=\frac{v}{2\pi}\log\left(z+\frac{r}{2N}\right)
  -\sigma\frac{(T_0L_0r)^{\frac{1}{2}}}{\pi (Nz)^{\frac{1}{2}}},
\end{align*}
and
\begin{align*}
  h(z)=h(z,r,t_g)=\frac{t_g}{2\pi}
  \log\left(1+\frac{r^2}{4N^2z^2}\left(1+\frac{r}{Nz}\right)^{-1}\right).
\end{align*}
It is easy to see that
$(z+\frac{r}{2N})^{u-1}$ is a flat function for $z$ and $r$.
Moreover, by the Taylor expansion and
$\frac{t_gr^2}{N^2z^2}\ll \frac{T_0R^2}{N^2}\ll \frac{T_0^3L^2C^4}{N^4}
\preccurlyeq \frac{T_0^3L^2}{T^4H^4}\preccurlyeq (TH)^{-\frac{1}{4}}$,
we see that $h(z)$ is a flat function for $z$ and $r$.
So $e(h(z))$ is also flat.
By applying Lemma \ref{lemma:stationary_phase BKY} with
\begin{align}\label{eqn:parameters}
  Y=S=\frac{(T_0LR)^{\frac{1}{2}}}{N^{\frac{1}{2}}}+|v|=\frac{T_0LC}{N}+|v|,
  \quad P=U=1, \quad X=T^\varepsilon,
\end{align}
one sees that it is negligible unless $|v|\asymp \frac{T_0LC}{N}$
and $\sgn(v)=\sgn(-\sigma)$.
In this case, by the Taylor expansion and
$\frac{vr}{N}\ll\frac{T_0^2L^2C^3}{N^3}\preccurlyeq \frac{T_0^2L^2}{T^3H^3}\preccurlyeq
T^{-\varepsilon}$,
we deduce that
\begin{align*}
  \phi_1(z)=\phi_1(z,r,v)=\frac{v}{2\pi}\left(\log\left(z+\frac{r}{2N}\right)-\log z\right)
\end{align*}
is a flat function for $z$, $r$ and $v$.
Let $V_0=\frac{T_0LC}{N}$.
Then, we smooth the parameter $v$ with $V(\frac{-\sigma v}{V_0})$
and move the weight function $V(\frac{r}{R})$ into the $z$-integral.
After separating the parameters
$z$, $r$ and $v$ in $(z+\frac{r}{2N})^{u-1}h(z)e(\phi_1(z))$
by using the technique in \S \ref{subsec:smooth weight functions},
$\tilde{G}_{r}(s)$ can be recast as
\begin{align*}
  4^{it_g}\int_0^\infty V\left(z,\frac{r}{R},\frac{v}{V_0}\right)e(\phi_2(z))\textup{d}z,
\end{align*}
where
\begin{align*}
  \phi_2(z)=\phi_2(z,r,v)
  =\frac{v}{2\pi}\log z - \sigma\frac{(T_0L_0r)^{\frac{1}{2}}}{\pi (Nz)^{\frac{1}{2}}},
\end{align*}
and $V(z,\frac{r}{R},\frac{v}{V_0})=V(z)V(\frac{r}{R})V(\frac{-\sigma v}{V_0})$.
By Lemma \ref{lem:main theorem in KPY} with
\begin{align*}
  X_1=1,\quad X_2=R,\quad X_3=Y=V_0,
\end{align*}
there is a stationary point
$z_0=\frac{T_0L_0r}{Nv^2}$
(actually $z_0^{\frac{1}{2}}=-\frac{\sigma}{v}(\frac{T_0L_0r}{N})^{\frac{1}{2}}$)
and $\tilde{G}_r(s)$ can be replaced by
$V_0^{-\frac{1}{2}}(\frac{T_0L_0r}{Nv^2})^{iv}e(\frac{v}{\pi})
  V(\frac{r}{R},\frac{v}{V_0})$,
where
$V(\frac{r}{R},\frac{v}{V_0})
=V(\frac{r}{R})V(\frac{-\sigma v}{V_0})$.
By absorbing the flat factor $(\frac{T_0L_0}{Nv^2})^{iv}e(\frac{v}{\pi})$
into $V(\frac{-\sigma v}{V_0})$,
we can rewrite $\tilde{G}_r(s)$ again as
\begin{align}\label{eqn: replace G}
  V_0^{-\frac{1}{2}}
  r^{iv}
  V\left(\frac{r}{R},\frac{v}{V_0}\right).
\end{align}
Note that \eqref{eqn: replace G} is not holomorphic with $s$,
so we will use this expression after moving the $v$-integration line.

We end this subsection by recall our parameters and conventions:
\begin{align}\label{recall parameters}
\begin{split}
  &t\in [T-HT^\varepsilon, T+HT^\varepsilon], \quad L=|t-t_g|, \quad H=o(L),
  \\
  &T_0=t+t_g\ll (TH)^{\frac{3}{4}-\varepsilon}, \quad L_0=|t-t_g|\asymp L,
   \quad N\preccurlyeq T_0L,
  \quad N_0=NT^{-\varepsilon},
  \\
  &C\preccurlyeq\frac{N}{TH}, \quad 0<r\asymp R=\frac{T_0LC^2}{N}\preccurlyeq N(TH)^{-\frac{1}{2}-\varepsilon},
  \quad
  |v|\asymp V_0=\frac{T_0LC}{N}.
\end{split}
\end{align}

\subsection{Treating $D^{\dag}_g(s,1,1,r)$}\label{subsec:deal hypergeometry}
Let $D^{\dag}_{g,F}(s,\nu_1,\nu_2,r)$ be the same as
$D_{g,F}(s,\nu_1,\nu_2,r)$ but restricted to $|m|\geq N_0$.
The purpose of this subsection is to see the relationship between
$D^{\dag}_g(s,1,1,r)$ and $D^{\dag}_{g,F}(s,1,1,r)$,
where $-\infty<\Re s<+\infty$ and the other related parameters
satisfy \eqref{recall parameters}.

We assume that the parameter $v$
satisfies $\kappa V_0\leq |v|\leq \kappa^{-1}V_0$, otherwise $\tilde{G}_r(s)$ is
negligible, where $\kappa$ is a sufficiently small positive number.
For the technical reason, we assume $\kappa^3V_0\leq |v|\leq \kappa^{-3}V_0$
in this subsection.
We first consider the hypergeometic function $F$ in $D_{g,F}$.
By \cite[9.100]{GraRyz}, we have
\begin{align}\label{eqn:def F}
  F\left(\frac{s}{2}+it_g,\frac{1}{2}+it_g,\frac{s+1}{2};
  \left(\frac{|m|-|n|}{|m|+|n|}\right)^2\right)
  =1+\sum_{\ell\geq1}\frac{(\frac{s}{2}+it_g)_\ell(\frac{1}{2}+it_g)_\ell}{(\frac{s+1}{2})_\ell\ell!}
  \left(\frac{|m|-|n|}{|m|+|n|}\right)^{2\ell},
\end{align}
where
\begin{align*}
  (a)_\ell=a(a+1)\cdots (a+\ell-1)
\end{align*}
is the Pochhammer symbol.
Note that we have
$
  |\frac{s}{2}+it_g+\ell|
  \leq 10(t_g+|v|+\left|\ell+\frac{u}{2}|\right),
$
$
  |\frac{1}{2}+it_g+\ell|
  \leq 10 (t_g+\ell),
$
and
$
  |\frac{s+1}{2}+\ell|
  \geq \frac{1}{10}(|v|+|\ell+\frac{u}{2}|),
$
uniformly in
\begin{align}\label{eqn:parameters 4.1}
  \kappa^3 V_0\leq |v|\leq \kappa^{-3}V_0, \quad \ell\geq 0, \quad u\in \mathbb{R}.
\end{align}

Thus, we can get
\begin{align*}
  \frac{(\frac{s}{2}+it_g)_\ell(\frac{1}{2}+it_g)_\ell}{(\frac{s+1}{2})_\ell\ell!}
  \leq \left(\frac{BT_0^2}{V_0}\right)^\ell,
\end{align*}
where $B$ is an absolutely positive large constant.
By using this, we rewrite \eqref{eqn:def F} as
\begin{align}\label{eqn: re F}
  F\left(\frac{s}{2}+it_g,\frac{1}{2}+it_g,\frac{s+1}{2};
  \left(\frac{|m|-|n|}{|m|+|n|}\right)^2\right)
  =\sum_{\ell\geq 0}a_\ell(s)w^\ell,
\end{align}
where $a_0(s)=1$, and
\begin{align*}
  a_\ell(s)=\frac{(\frac{s}{2}+it_g)_\ell(\frac{1}{2}+it_g)_\ell}{(\frac{s+1}{2})_\ell\ell!}
  \left(\frac{2BT_0^2}{V_0}\right)^{-\ell}
\end{align*}
for $\ell\geq 1$, and
\begin{align*}
  w=\frac{2BT_0^2}{V_0}\left(\frac{|m|-|n|}{|m|+|n|}\right)^{2}.
\end{align*}
Note that we have
$|a_\ell(s)|\leq \frac{1}{2^\ell}$ and
\begin{align*}
  w \ll \frac{T_0^2R^2}{V_0N_0^2}=T^{2\varepsilon}\frac{T_0^3LC^3}{N^3}
  \preccurlyeq T^{-\varepsilon}.
\end{align*}
Now let $b_0(s)=1$
and
\begin{align}\label{eqn:def b}
  b_\ell(s)=-\sum_{\substack{\ell_1+\ell_2=\ell\\ 0\leq\ell_1<\ell}}
  b_{\ell_1}(s)a_{\ell_2}(s+2\ell_1),
\end{align}
for $\ell\geq 1$.
Clearly, in the region $\kappa^3V_0\leq |v|\leq \kappa^{-3}V_0$,
$b_\ell(s)$ is holomorphic and satisfies the upper bound
$|b_\ell(s)|\leq 1$.
We claim that one can choose a sufficiently large positive $A$ which depends only
on $\varepsilon$, such that
we can approximate $1$ by
\begin{align}\label{eqn:express 1 by F}
  \sum_{0\leq \ell\leq A}b_\ell(s) F\left(\frac{s+2\ell}{2}+it_g,\frac{1}{2}+it_g,\frac{s+2\ell+1}{2};
  \left(\frac{|m|-|n|}{|m|+|n|}\right)^2\right)w^\ell
\end{align}
at the cost of a negligible error term.
Actually, we have
\begin{multline*}
  \left|\sum_{0\leq \ell_1\leq A}b_{\ell_1}(s)(\sum_{0\leq \ell_2\leq A}a_{\ell_2}(s+2\ell_1)w^{\ell_2})w^{\ell_1}
  \right.
  \\
  \left.-\sum_{0\leq \ell_1\leq A}b_{\ell_1}(s) F\left(\frac{s+2\ell_1}{2}+it_g,\frac{1}{2}+it_g,\frac{s+2\ell_1+1}{2};
  \left(\frac{|m|-|n|}{|m|+|n|}\right)^2\right)w^{\ell_1}\right|
  \\
  \ll \sum_{0\leq \ell_1\leq A}|b_{\ell_1}(s)|w^{\ell_1}
  \sum_{\ell_2\geq A}  |a_{\ell_2}(s+2\ell_1)|w^{\ell_2}
  \ll \sum_{0\leq \ell_1\leq A}|b_{\ell_1}(s)|w^{\ell_1}T^{-A\varepsilon}
  \ll T^{-A\varepsilon},
\end{multline*}
and
\begin{align*}
  &\sum_{0\leq \ell_1\leq A}b_{\ell_1}(s)(\sum_{0\leq \ell_2\leq A}a_{\ell_2}(s+2\ell_1)w^{\ell_2})w^{\ell_1}
  \\
  =&\sum_{0\leq \ell\leq A}w^\ell\left(
  \sum_{\substack{\ell_1+\ell_2=\ell\\ \ell_1,\,\ell_2\geq 0}}b_{\ell_1}(s)a_{\ell_2}(s+2\ell_1)\right)
  +\sum_{A<\ell\leq 2A}w^\ell\left(
  \sum_{\substack{\ell_1+\ell_2=\ell\\ 0\leq \ell_1,\,\ell_2\leq A}}b_{\ell_1}(s)a_{\ell_2}(s+2\ell_1)\right)
  \\
  =&1+O(T^{-A\varepsilon}).
\end{align*}
So the claim follows by letting $A=\frac{1000}{\varepsilon}$, say.
Consequently, by noting that $(|m|-|n|)^2=r^2$ when $|m|>N_0$
and $n-m=r$,
we insert \eqref{eqn:express 1 by F} into $D_g^{\dag}(s,1,1,r)$
and deduce that
\begin{align}\label{eqn:express D dag}
  D_g^{\dag}(s,1,1,r)=\sum_{0\leq \ell\leq A}b_{\ell}(s)\left(\frac{2BT_0^2r^2}{V_0}\right)^{\ell}D_{g,F}^{\dag}(s+2\ell,1,1,r)
  +O(T^{-1000}).
\end{align}


Denote by
\begin{align}\label{eqn:integration line gamma}
  \gamma(u)=\{s=u+iv|\kappa V_0\leq |v|\leq \kappa^{-1} V_0\}.
\end{align}
Then, by \eqref{eqn:express D dag} and \eqref{eqn:express P by D}, we get
\begin{align}\label{eqn:express P by D 3}
  P(r,N)=\sum_{\ell\leq A}\left(\frac{BT_0^2r^2}{2V_0N^2}\right)^{\ell}
  (P_{1,\ell}(r,N)-P_{2,\ell}(r,N))
  +O(T^{-1000}),
\end{align}
where
\begin{align}\label{P1}
  P_{1,\ell}(r,N)=\frac{1}{4\pi i}\int_{\gamma(2)}b_{\ell}(s)
  (2N)^{s+2\ell}\tilde{G}_{r}(s)D_{g,F}(s+2\ell,1,1,r)\textup{d}s
\end{align}
and
\begin{align}\label{P2}
  P_{2,\ell}(r,N)=\frac{1}{4\pi i}\int_{\gamma(2)}b_{\ell}(s)
  (2N)^{s+2\ell}\tilde{G}_{r}(s)D_{g,F}^{\ddag}(s+2\ell,1,1,r)\textup{d}s,
\end{align}
with $D_{g,F}^{\ddag}(s,\nu_1,\nu_2,r)$ being the same as $D_{g,F}(s,\nu_1,\nu_2,r)$
but restricted as $|m|\leq N_0$.
\subsection{The contribution of $P_{2,\ell}(r,N)$}
\label{subsec:treat P2}
For our convenience, we denote by $z=1-\frac{n^2}{m^2}$ in this
subsection.
By the definition of $D^{\ddag}_{g,F}(s,\nu_1,\nu_2,r)$ and \eqref{eqn:eqn of F}, we get
\begin{multline}\label{eqn:re P2}
  P_{2,\ell}(r,N)= \sum_{\substack{n-m=r\\ |m|\leq N_0}}\lambda_g(m)\lambda_g(n)
  \left(\frac{|n|}{4|m|}\right)^{it_g}
  \\
  \cdot \frac{1}{4\pi i}\int_{\gamma(2\ell+2)}
  b_\ell(s-2\ell)\tilde{G}_{r}(s-2\ell)\left(\frac{N}{|m|}\right)^{s}
  F\left(\frac{s}{2}+it_g,\frac{s}{2},s;z\right)ds,
\end{multline}
where we have replaced $s+2\ell$ by $s$.
Note that we also have $|n|\ll N_0$.
We will show the $s$-integral in \eqref{eqn:re P2} is negligible.
By (see \cite[9.111]{GraRyz})
\begin{align}\label{eqn:pro of F}
  F(\alpha,\beta,\gamma;z)=\frac{1}{B(\beta, \gamma-\beta)}
  \int_0^1 t^{\beta-1}(1-t)^{\gamma-\beta-1}(1-tz)^{-\alpha}dt, \quad \Re\gamma>\Re\beta>0,
\end{align}
and (see \cite[8.384-1]{GraRyz})
\begin{align}\label{eqn:pro of B}
  B(x,y)=\frac{\Gamma(x)\Gamma(y)}{\Gamma(x+y)}, \quad \Re x>0, \Re y>0,
\end{align}
one obtains
\begin{align}\label{eqn:t integral in F}
  F\left(\frac{s}{2}+it_g,\frac{s}{2},s;z\right)
  =\frac{\Gamma(s)}{\Gamma^2(\frac{s}{2})}
  \int_0^1 t^{\frac{s}{2}-1}(1-t)^{\frac{s}{2}-1}
  (1-tz)^{-\frac{s}{2}-it_g}dt.
\end{align}
The strategy is that, we would like to remove the $s$-integral to far left,
so that one might see the contribution is negligible from the factor $\frac{N}{|m|}$
and $|m|\leq N_0$
(note that there is no pole in this process since $|\Im s|\asymp V_0$).
The problem is that the $t$-integral in \eqref{eqn:t integral in F} is not absolutely convergent
at $t=0$, $1$ when $\Re s\leq 0$, so that the integration line can only be moved to $\Re s=\varepsilon$.
Nevertheless, we can move the line more left after
partial integration with respect to $t$ repeatedly.

It is obvious that, for $0\leq t\leq 1$, one has
$
  \min\{1,\frac{n^2}{m^2}\}\leq 1-tz
  \leq \max\{1,\frac{n^2}{m^2}\},
$
which implies that
\begin{align}\label{eqn:bound 1-tz}
  (1-tz)^s\ll 1+\left(\frac{n^2}{m^2}\right)^{\Re s}.
\end{align}
Now we consider
\begin{align}\label{eqn: t-integral 1}
  \int_{\frac{1}{2}}^1t^{\frac{s}{2}-1}(1-t)^{\frac{s}{2}-1}
  (1-tz)^{-\frac{s}{2}-it_g}dt,
\end{align}
which only has the convergence problem at $t=1$.
By partial integration with respect to $(1-t)^{\frac{s}{2}-1}$, it becomes
(note that $\sigma_3=\sgn(n-m)$)
\begin{align}\label{eqn:partial once}
  &\frac{2^{2-s}}{s}\left(1-\frac{z}{2}\right)^{-\frac{s}{2}-it_g}\nonumber
  \\
  &+\frac{2}{s}\left(\frac{s}{2}-1\right)\int_{\frac{1}{2}}^1 (1-t)^{\frac{s}{2}}
  t^{\frac{s}{2}-2}
  (1-tz)^{-\frac{s}{2}-it_g}dt\nonumber
  \\
  &-\frac{2}{s}\left(\frac{s}{2}+it_g\right)\left(\frac{r(m+n)}{m^2}\right)
  \int_{\frac{1}{2}}^1 (1-t)^{\frac{s}{2}}t^{\frac{s}{2}-1}
  (1-tz)^{-\frac{s}{2}-it_g-1}dt.
\end{align}
For the first term of \eqref{eqn:partial once}, the corresponding $s$-integral
in \eqref{eqn:re P2} can be treated by moving the integration line to far left.
Actually, by \eqref{eqn:bound 1-tz}, we have, for $\Re s=-2A$,
\begin{align*}
  \left(\frac{N}{|m|}\right)^s\left(1-\frac{z}{2}\right)^{-\frac{s}{2}-it_g}
  \ll \left(\frac{|m|}{N}\right)^{2A}
  +\left(\frac{|m|}{N}\right)^{2A}\left(\frac{n^2}{m^2}\right)^A
  \ll T^{-A\varepsilon}.
\end{align*}
Thus the corresponding $s$-integral is negligible
by taking $A$ large enough.

For the rest two terms of \eqref{eqn:partial once},
the $t$-integral is absolutely convergent at $1$ when $\Re s\geq -2+\varepsilon$,
so the integration line of resulting $s$-integral can be moved to $\Re s=-2+\varepsilon$.
We can repeat the partial integration $A$ times.
Those terms like the first term of \eqref{eqn:partial once} can be treated similarly:
move the integration line to far left, and the contribution can be omitted.
For the rest terms, we only deal with the two typical terms,
and the corresponding $s$-integral is bounded by
\begin{align*}
  \int_{\frac{1}{2}}^1 \left|\int_{\gamma(2\ell+2)}
  b_\ell(s-2\ell)\tilde{G}_{r}(s-2\ell)\left(\frac{N}{|m|}\right)^{s}
  \frac{\Gamma(s)}{\Gamma^2(\frac{s}{2})}(|I_1|+|I_2|)
  ds\right|dt,
\end{align*}
where
\begin{align*}
  I_1=\prod_{j=1}^A\left(\frac{s}{2}+j-1\right)^{-1}\left(\frac{s}{2}-j\right)
  (1-t)^{\frac{s}{2}+A-1} t^{\frac{s}{2}-A-1}
  (1-tz)^{-\frac{s}{2}-it_g},
\end{align*}
and
\begin{align*}
  I_2= \prod_{j=1}^A\left(\frac{s}{2}+j-1\right)^{-1}\left(\frac{s}{2}+it_g+j-1\right)
  \left(\frac{r(m+n)}{m^2}\right)^A
  \\
  \cdot (1-t)^{\frac{s}{2}+A-1} t^{\frac{s}{2}-1}
  (1-tz)^{-\frac{s}{2}-it_g-A}.
\end{align*}
Moving the $s$-integral to $\Re s=-2A+2\varepsilon$, we have
\begin{align*}
  \left(\frac{N}{|m|}\right)^s|I_1|\preccurlyeq
  (1-t)^{\varepsilon-1}\left(\left(\frac{|m|}{N}\right)^{2A}
  +\left(\frac{|m|}{N}\right)^{2A}\left(\frac{n^2}{m^2}\right)^A\right)
  \preccurlyeq (1-t)^{\varepsilon-1}T^{-A\varepsilon},
\end{align*}
and
\begin{align*}
  \left(\frac{N}{|m|}\right)^s|I_2|\preccurlyeq
   (1-t)^{\varepsilon-1}
  \left(\frac{|m|}{N}\right)^{2A}\left(\frac{T_0RN_0}{V_0m^2}\right)^{A}
  \preccurlyeq (1-t)^{\varepsilon-1}(TH)^{-(\frac{1}{4}+\varepsilon)A},
\end{align*}
which implies that the contribution is negligible.

For $0\leq t\leq \frac{1}{2}$, by partial integration with respect to $t^{\frac{s}{2}-1}$
and the similar treatment, the contribution of
\begin{align*}
  \int_{0}^{\frac{1}{2}}t^{\frac{s}{2}-1}(1-t)^{\frac{s}{2}-1}
  (1-tz)^{-\frac{s}{2}-it_g}dt
\end{align*}
can also be omitted.


\subsection{The contribution of $P_{1,\ell}(r,N)$}
By \eqref{eqn:decomposition of D},
\eqref{P1}
and replace $s+2\ell$ with $s$, we get
\begin{align*}
  P_{1,\ell}(r,N)=P_{1,\ell,d}(r,N)+P_{1,\ell,E}(r,N),
\end{align*}
where
\begin{align*}
  P_{1,\ell,d}(r,N)=\frac{1}{4\pi i}\int_{\gamma(2\ell+2)}b_{\ell}(s-2\ell)
  (2N)^{s}\tilde{G}_{r}(s-2\ell)D_{g,F,d}(s,1,1,r)ds,
\end{align*}
and
\begin{align*}
  P_{1,\ell,E}(r,N)=\frac{1}{4\pi i}\int_{\gamma(2\ell+2)}b_{\ell}(s-2\ell)
  (2N)^{s}\tilde{G}_{r}(s-2\ell)D_{g,F,E}(s,1,1,r)ds.
\end{align*}
By Lemma \ref{lem:inner product exchange}, we only need to consider
$\phi\in \mathcal{B}^*(q_1)$ and $E_{\mathfrak{a}}=E_\infty$ in \eqref{eqn:Dd}
and \eqref{eqn:DE}, respectively, where $q_1|q$.
By Lemma \ref{lem:Watson-Ichino q}, we have
\begin{align*}
  |\left<\phi,|g|^2\right>_q|^2\preccurlyeq
  \left|\frac{\Lambda(\frac{1}{2},\phi)\Lambda(\frac{1}{2},\phi\times \textup{ad}\,g)}{\Lambda(1,\textup{ad}\,g)^2\Lambda(1,\textup{sym}^2\,\phi)}\right|.
\end{align*}
As we said in \S 2.2, we can assume $\phi$ is even.
By \eqref{eqn:ad g}-\eqref{eqn:zeta value at 1}
and the Stirling formula, we get
\begin{align}\label{eqn:triple product}
  \left<\phi,|g|^2\right>_q\preccurlyeq
  \frac{e^{-\frac{\pi}{2} \Omega(t_\phi,t_g)}
  \sqrt{L(\frac{1}{2},\phi)L(\frac{1}{2},\phi\times \textup{ad}g)}}
  {(1+t_\phi)^{\frac{1}{2}}
  (1+|2t_g+t_\phi|)^{\frac{1}{4}}(1+|2t_g-t_\phi|)^{\frac{1}{4}}},
\end{align}
where
\begin{align}\label{eqn:define omega}
  \Omega(t_1,t_2)=\begin{cases}
        0,\quad & 0\leq |t_1|\leq 2|t_2|,
        \\\noalign{\vskip 1mm}
        |t_1|-2|t_2|, \quad &0\leq 2|t_2|\leq |t_1|.
     \end{cases}
\end{align}
Moreover, by the Stirling formula again together with \eqref{eqn:Ag}, \eqref{eqn:B} and
\eqref{eqn:Fc3},
one gets (for $\mu\in\mathbb{R}\cup i(-\frac{1}{2},\frac{1}{2})$)
\begin{align}\label{eqn:bound A}
  A_{g}(s)&\preccurlyeq e^{-\frac{\pi}{2}(2t_g-|t_g+\frac{v}{2}|-|t_g-\frac{v}{2}|)}
  (1+|2t_g+v|)^{\frac{1-u}{2}}(1+|2t_g-v|)^{\frac{1-u}{2}}
  V_0^{1/2}\nonumber
  \\
  &\preccurlyeq e^{\frac{\pi}{2}\Omega(v,t_g)}
  (1+|2t_g+v|)^{\frac{1-u}{2}}(1+|2t_g-v|)^{\frac{1-u}{2}}V_0^{1/2},
\end{align}
and
\begin{align}\label{eqn:bound B}
  B(s,\mu)&\ll e^{\frac{\pi}{2}(|\mu|-\frac{|\mu+v|-|\mu-v|}{2})}
  (1+|\mu+v|)^{\frac{u}{2}-\frac{3}{4}}
  (1+|\mu-v|)^{\frac{u}{2}-\frac{3}{4}}\nonumber
  \\
  &\ll e^{-\frac{\pi}{2}\Omega(v,\frac{\mu}{2})}
  (1+|\mu+v|)^{\frac{u}{2}-\frac{3}{4}}
  (1+|\mu-v|)^{\frac{u}{2}-\frac{3}{4}}.
\end{align}
By some simple calculations, we get
\begin{align}\label{eqn:Omega}
  e^{\frac{\pi}{2}\Omega(v,t_g)}e^{-\frac{\pi}{2}\Omega(v,\frac{t_\phi}{2})}
  e^{-\frac{\pi}{2} \Omega(t_\phi,t_g)}
  \ll \begin{cases}
        \min\{1, e^{-\frac{\pi}{2}(|t_\phi|-2t_g)}\}
        \quad & |v|\leq 2t_g,
        \\\noalign{\vskip 1mm}
        e^{-\frac{\pi}{2}(|t_\phi|-|v|)} \quad & 2t_g<|v|\leq t_\phi,
        \\\noalign{\vskip 1mm}
        e^{-\frac{\pi}{2}(2t_g-|t_\phi|)} \quad & t_\phi<2t_g<|v|,
        \\\noalign{\vskip 1mm}
        1 \quad & 2t_g\leq t_\phi<|v|,
        \\\noalign{\vskip 1mm}
        \max\{e^{-\pi t_g}, e^{-\frac{\pi}{2}|v|}\}
        \quad & t_\phi\in i(-\frac{1}{2},\frac{1}{2}).
     \end{cases}
\end{align}
So, by noting that $\kappa V_0\leq |v|\leq \kappa^{-1}V_0$, we may first truncate $t_\phi$ at
$t_\phi\in i(-\frac{1}{2},\frac{1}{2})$ or
$t_\phi\leq \max\{2t_g+t_g^{1-\varepsilon},\kappa^{-2}V_0\}$
for some positive number $\varepsilon$, since otherwise, we have an exponential decay.
Note that, by \eqref{eqn:Dd}, the $s$-integral in $D_{g,F,d}$ is
\begin{align}\label{eqn:s integal in DgFd}
  \int_{\gamma(2\ell+2)}b_\ell(s-2\ell)(2N)^s\tilde{G}_r(s-2\ell)
  \frac{A_{g}(s)B(s,t_j)}{r^{s-1}}\left<\phi,|g|^2\right>_qds.
\end{align}
Hence, if $t_\phi\in (0,\kappa^2V_0)\cup (\kappa^{-2}V_0,2t_g+t_g^{1-\delta})
\cup i(-\frac{1}{2},\frac{1}{2})$,
then we move the integration line
to far left at the cost of a negligible error term. In fact, recall $b_\ell(s)$ is holomorphic in
$\kappa^3 V_0\leq |v|\leq \kappa^{-3} V_0$,
no pole is encountered during this shifting and
$\tilde{G}_{r}(s)$ is arbitrary small on the horizontal line segments.
For $\Re s=-A$, we have, by using \eqref{recall parameters},
\eqref{eqn:triple product},
\eqref{eqn:bound A}, \eqref{eqn:bound B} and \eqref{eqn:Omega},
\begin{align*}
  e^{-\frac{\pi}{2} \Omega(t_{\phi},t_g)}N^s \tilde{G}_r(s-2\ell)
  \frac{A_{g}(s)B(s,t_{\phi})}{r^{s-1}}
  \ll T_0RV_0^{-\frac{3}{2}}\left(\frac{T_0R}{NV_0}\right)^A
  \preccurlyeq T_0RV_0^{-\frac{3}{2}}(TH)^{(\frac{1}{4}-\varepsilon)A}.
\end{align*}
So one sees that \eqref{eqn:s integal in DgFd} is negligible.
This allows us to finally truncate $t_\phi$ at $\kappa^2V_0\leq t_\phi\leq \kappa^{-2}V_0$.
Similarly, we can truncate $\tau$ at $\kappa^2V_0\leq |\tau|\leq \kappa^{-2}V_0$.

\subsubsection{The contribution of discrete spectrum}\label{subsubsec:contribution d}
By the above argument together with
\eqref{eqn:express P by D 3}, \eqref{P1}, \eqref{rere goal Maass},
\eqref{eqn:replace the c sum}, \eqref{eqn:decomposition of D} and
\eqref{eqn:Dd}, we need to get the following estimate
\begin{multline}\label{eqn:goal Pd}
  \frac{TH}{C^{3/2}}\left(\frac{N}{T_0L}\right)^{\frac{1}{2}}
  \sum_{\ell\leq A}\left(\frac{T_0^2R^2}{V_0N^2}\right)^{\ell}
  \sum_{1\leq k\ll \log T}\sum_{r\geq 1}\mathcal{S}_{r,k}
  \sum_{\substack{\phi\in \mathcal{B}^*(q_1)
  \\ \kappa^2V_0 \leq t_\phi\leq \kappa^{-2}V_0}}
  \frac{\overline{\rho}_\phi(r)}{\cosh(\frac{\pi t_\phi}{2})}\left<\phi,|g|^2\right>_q
  \\
  \cdot
  \int_{\gamma(2\ell+2)}b_\ell(s-2\ell)(2N)^s\tilde{G}_r(s-2\ell)
  \frac{A_{g}(s)B(s,t_\phi)}{r^{s-1}}ds
  \preccurlyeq NT_0L^{\frac{1}{2}}||g||^2_4+NTH,
\end{multline}
where
\begin{align*}
  \mathcal{S}_{r,k}=V\left(\frac{r}{R}\right)G_{\chi\psi}(r,[q,D]^k)
  \sum_{d|r}d\sum_{\substack{c_k\geq 1\\d|c_k}}\mu\left(\frac{c_k}{d}\right)
  V\left(\frac{c_k}{C_k}\right)\chi\psi(c_k).
\end{align*}

We may drop the $k$- and $\ell$-sums, since $k\ll \log T$
and $\frac{T_0^2R^2}{V_0N^2}\preccurlyeq (TH)^{-\varepsilon} $.
Now we shift the contour to $\Re s=-\frac{1}{2}$.
It is easy to see that it has at most one simple pole at
$s=\frac{1}{2}-i\sigma t_\phi$ (recall $\textup{sgn}(v)=\textup{sgn}(-\sigma)$).
Note that this requires in addition to the vanishing of $\Gamma^{-1}(\frac{s}{2})$
at $s=0$.
We may also assume that no $t_\phi$ is equal to $\kappa V_0$ or $\kappa^{-1}V_0$ by a small perturbation
(of magnitude $\varepsilon$).
By \eqref{eqn: replace G}, \eqref{eqn:Ag} and \eqref{eqn:B},
the residue can be replaced with
\begin{multline}\label{eqn:replace residue}
  b_\ell\left(\frac{1}{2}-2\ell-i\sigma t_\phi\right)
  (2N)^{\frac{1}{2}-i\sigma t_\phi}
  r^{\frac{1}{2}}V_0^{-\frac{1}{2}}\left(\frac{T_0L_0}{Nt_\phi^2}\right)^{-i\sigma t_\phi}
  \\
  \cdot A_g\left(\frac{1}{2}-i\sigma t_\phi\right)
  \Gamma(-i\sigma t_\phi)\cosh\left(\frac{\pi t_\phi}{2}\right),
\end{multline}
which can be bounded by
\begin{align*}
  e^{\frac{\pi}{2}\Omega(t_\phi,t_g)}\left(\frac{NR}{V_0}\right)^{\frac{1}{2}}
  (t_g+V_0)^{\frac{1}{4}}(1+|2t_g-t_\phi|)^{\frac{1}{4}}.
\end{align*}
Note that $r$ and $t_\phi$ have been separated.
Hence, to compute its contribution to \eqref{eqn:goal Pd},
we are led to consider
\begin{multline}\label{eqn:goal Pd re}
  \frac{TH}{C^{3/2}}\left(\frac{N}{T_0L}\right)^{\frac{1}{2}}
  \left(\frac{NR}{V_0}\right)^{\frac{1}{2}}
  (t_g+V_0)^{\frac{1}{4}}
  \\
  \cdot \sum_{\substack{\phi\in \mathcal{B}^*(q_1)
  \\ \kappa^2V_0 \leq t_\phi\leq \kappa^{-2}V_0}}
  \left|e^{\frac{\pi}{2}\Omega(t_\phi,t_g)}(1+|2t_g-t_\phi|)^{\frac{1}{4}}
  \left<\phi,|g|^2\right>_q\right|
  \left|\sum_{r\asymp R}
  \frac{\bar{\rho}_\phi(r)\mathcal{S}_{r,k}}{\cosh(\frac{\pi t_\phi}{2})}\right|.
\end{multline}
By the Cauchy--Schwarz inequality, the second line of \eqref{eqn:goal Pd re}
is bounded by $(\mathfrak{C}_1\mathfrak{C}_2)^{1/2}$,
where
\begin{align*}
  \mathfrak{C}_1
  =\sum_{\substack{\phi\in \mathcal{B}^*(q_1)
  \\ \kappa^2V_0 \leq t_\phi\leq \kappa^{-2}V_0}}
  \frac{1}{\cosh \pi t_\phi}
  \left|\sum_{r\asymp R}\bar{\rho}_\phi(r)\mathcal{S}_{r,k}\right|^2,
\end{align*}
and
\begin{align*}
  \mathfrak{C}_2=\sum_{\substack{\phi\in \mathcal{B}^*(q_1)
  \\ \kappa^2V_0 \leq t_\phi\leq \kappa^{-2}V_0}}
  e^{\Omega(t_\phi,t_g)}(1+|2t_g-t_\phi|)^{\frac{1}{2}}
  |\left<\phi,|g|^2\right>_q|^2.
\end{align*}
By the large sieve inequality (Lemma \ref{lem:large sieve})
and $\mathcal{S}_{r,k}\preccurlyeq [q,D]^kC_k=C$, we have
\begin{align*}
  \mathfrak{C}_1\preccurlyeq C^2(V_0^2+R)\preccurlyeq C^2V_0^2.
\end{align*}
So we obtain
\begin{align*}
  \eqref{eqn:goal Pd re}\ll &
  \frac{TH}{C^{3/2}}\left(\frac{N}{T_0L}\right)^{\frac{1}{2}}
  \left(\frac{NR}{V_0}\right)^{\frac{1}{2}}
  (t_g+V_0)^{\frac{1}{4}}
  CV_0
  \mathfrak{C}_2^{\frac{1}{2}}
  \\
  \ll & THC(T_0L)^{\frac{1}{2}}(t_g+V_0)^{\frac{1}{4}}
  \mathfrak{C}_2^{\frac{1}{2}}.
\end{align*}
Now we turn to $\mathfrak{C}_2$.
If $t_g\gg \frac{T^{1+\varepsilon}}{H}$, then, by noting that
\begin{align*}
  V_0\preccurlyeq\frac{T_0L}{TH}\preccurlyeq\frac{t_g^2+T^2}{TH}
  \preccurlyeq \frac{t_g}{(TH)^{\frac{1}{4}}}+\frac{T}{H},
\end{align*}
one has
$t_g\gg V_0^{1+\varepsilon}$.
We thus obtain
\begin{align*}
  \mathfrak{C}_2\ll T_0^{\frac{1}{2}}||g||_4^4.
\end{align*}
It follows that
\begin{align*}
  \eqref{eqn:goal Pd re}\ll THC(T_0L)^{\frac{1}{2}}(t_g+V_0)^{\frac{1}{4}}
  (T_0^{\frac{1}{2}}||g||_4^4)^{\frac{1}{2}}
  \preccurlyeq NT_0L^{\frac{1}{2}}||g||_4^2.
\end{align*}

If $t_g\preccurlyeq \frac{T}{H}=o(T_0)$,
then we have $T_0\asymp T\asymp L$ and
$V_0\preccurlyeq\frac{T_0L}{TH}\preccurlyeq\frac{T}{H}$.
Therefore, by \eqref{eqn:triple product} and  the Cauchy-Schwarz inequality,
$\mathfrak{C}_2$ can be bounded by
\begin{align}\label{eqn:reduce to triple product}
  \preccurlyeq V_0^{-1}(t_g+V_0)^{-\frac{1}{2}}
  \mathfrak{C}_3^{\frac{1}{2}}\mathfrak{C}_4^{\frac{1}{2}},
\end{align}
where
\begin{align*}
  \mathfrak{C}_3=\sum_{\substack{\phi\in \mathcal{B}^*(q_1)
  \\ t_\phi\asymp V_0}}
  \left|L\left(\frac{1}{2},\phi\times \textup{ad}\,g\right)\right|^2,
\end{align*}
and
\begin{align*}
  \mathfrak{C}_4=\sum_{\substack{\phi\in \mathcal{B}^*(q_1)
  \\ t_\phi\asymp V_0}}
  \left|L\left(\frac{1}{2},\phi\right)\right|^2.
\end{align*}
To bound this, we quote the following result.
\begin{lemma}\label{lem:HK22 lemma 6.1}
\footnote{In \cite[Proposition 6.1]{HumKhan22},
the authors got this result when $q=1$.
It is clear to see that
the tools there
(approximate functional equation, the large sieve inequality,
Li's result \cite[Theorem 2]{Li10}, etc.) work well in our case.}
\cite[Proposition 6.1]{HumKhan22}
  For $M\geq U\geq 1$, we have
  \begin{align*}
   \sum_{\substack{\phi\in \mathcal{B}^*(q_1)
  \\ M\leq t_\phi\leq 2M}}
  \frac{|L(\frac{1}{2},\phi\times \textup{ad}\,g)|^2}
  {L(1,\textup{ad}\,g)}+
  \underset{M\leq |\tau|\leq 2M}{\int}
  \left|\frac{L(\frac{1}{2}+it,\textup{ad}\,g)}{\zeta(1+2it)}\right|^2dt
  \ll
  \begin{cases}
    t_g^{2+\varepsilon}M, & \mbox{if $M\leq 2t_g$,}  \\
    M^{3+\varepsilon}, & \mbox{if $M\geq 2t_g$},
  \end{cases}
  \end{align*}
  and
   \begin{align*}
   \sum_{\substack{\phi\in \mathcal{B}^*(q_1)
  \\ M-U\leq t_\phi\leq M+U}}
  \frac{|L(\frac{1}{2},\phi)|^2}
  {L(1,\textup{ad}\,g)}+
  \underset{M-U\leq |\tau|\leq M+U}{\int}
  \left|\frac{\zeta(\frac{1}{2}+it)^2}{\zeta(1+2it)}\right|^2dt
  \ll M^{1+\varepsilon}U.
  \end{align*}
\end{lemma}
Now we come back to estimate $\mathfrak{C}_2$.
By Lemma \ref{lem:HK22 lemma 6.1} and \eqref{eqn:reduce to triple product},
we obtain
\begin{align*}
  \mathfrak{C}_2\preccurlyeq
  V_0^{-1}(t_g+V_0)^{-\frac{1}{2}}
  (V_0^3(t_g^2+V_0^2))^{\frac{1}{2}}
  \preccurlyeq \frac{T}{H}.
\end{align*}
Hence, we deduce that
\begin{align}\label{eqn:tg small residue}
  \eqref{eqn:goal Pd re}
  &\preccurlyeq \frac{TH}{C^{3/2}}\left(\frac{N}{T_0L}\right)^{\frac{1}{2}}
  \left(\frac{NR}{V_0}\right)^{\frac{1}{2}}
  (t_g+V_0)^{\frac{1}{4}}
  CV_0\left(\frac{T}{H}\right)^{\frac{1}{2}}\nonumber
  \\
  &\preccurlyeq NT^{\frac{3}{2}}\left(\frac{T}{H^3}\right)^\frac{1}{4}\nonumber
  \\
  &\preccurlyeq NT^{\frac{3}{2}},
\end{align}
which is acceptable.

To get \eqref{eqn:goal Pd},
we still need to evaluate the contribution coming from the
$s$-integral when $\Re s=-\frac{1}{2}$.
By \eqref{eqn: replace G}, we may replace $\frac{\tilde{G}_r(-1/2-2\ell)}{r^{s-1}}$
with
$R^{\frac{3}{2}}V_0^{-\frac{1}{2}}(\frac{T_0L_0}{Nv^2})^{iv}$,
which is bounded by $R^{\frac{3}{2}}V_0^{-\frac{1}{2}}$.
We emphasize again that $r$ and $v$ have been separated.
Moreover, we apply \eqref{eqn:bound A}-\eqref{eqn:bound B} and get
\begin{align*}
  A_g\left(-\frac{1}{2}+iv\right)B\left(-\frac{1}{2}+iv,t_\phi\right)
  \preccurlyeq
  e^{\frac{\pi}{2}(\Omega(v,t_g)-\Omega(v,\frac{t_\phi}{2}))}
  (t_g+V_0)^{\frac{3}{2}}V_0^{-\frac{1}{2}}
  (1+|t_\phi-|v||)^{-1}.
\end{align*}
So, by the Cauchy--Schwarz inequality,
it is enough to bound
\begin{align}\label{eqn:goal Pd -1/2}
  \frac{TH}{C^{3/2}}\left(\frac{N}{T_0L}\right)^{\frac{1}{2}}
  R^{\frac{3}{2}}(t_g+V_0)^{\frac{3}{2}}N^{-\frac{1}{2}}V_0^{-1}
  (\mathfrak{J}_1\mathfrak{J}_2)^{\frac{1}{2}},
\end{align}
by $\preccurlyeq NT_0L^{\frac{1}{2}}||g||_4^2+NTH$, where
\begin{align*}
  \mathfrak{J}_1=
  \sum_{t_\phi\asymp V_0}
  \frac{1}{\cosh \pi t_\phi}
  \left|\sum_{r\asymp R}
  \bar{\rho}_\phi(r)
  \mathcal{S}_{r,k}\right|^2
  \int_{\gamma(-\frac{1}{2})}(1+|t_\phi-|v||)^{-1}|ds|,
\end{align*}
and
\begin{align*}
  \mathfrak{J}_2=
  \sum_{t_\phi\asymp V_0}
  \left|\left<\phi,|g|^2\right>_q\right|^2
  \int_{\gamma(-\frac{1}{2})}
    e^{\frac{\pi}{2}(\Omega(v,t_g)-\Omega(v,\frac{t_\phi}{2}))}
    (1+|t_\phi-|v||)^{-1}|ds|.
\end{align*}
For $\mathfrak{J}_1$, we drop the $s$-integral
and use the large sieve inequality (Lemma \ref{lem:large sieve}),
obtaining
\begin{align*}
  \mathfrak{J}_1\preccurlyeq (V_0^2+R)C^2\preccurlyeq C^2V_0^2.
\end{align*}
If $t_g\gg \frac{T^{1+\varepsilon}}{H}$,
one has (note that we have got $t_g\gg V_0^{1+\varepsilon}$)
\begin{align*}
  e^{\frac{\pi}{2}(\Omega(v,t_g)-\Omega(v,\frac{t_\phi}{2}))}\ll 1,
\end{align*}
and hence
\begin{align*}
  \mathfrak{J}_2\ll ||g||^4_4.
\end{align*}
It follows that
\begin{align*}
  \eqref{eqn:goal Pd -1/2}
  &\preccurlyeq
  \frac{TH}{C^{3/2}}\left(\frac{N}{T_0L}\right)^{\frac{1}{2}}
  R^{\frac{3}{2}}(t_g+V_0)^{\frac{3}{2}}N^{-\frac{1}{2}}V_0^{-1}
  CV_0||g||_4^2
  \\
  &\preccurlyeq
  NT_0L^{\frac{1}{2}}\frac{T_0^{\frac{3}{2}}L^{\frac{1}{2}}}{(TH)^{\frac{3}{2}}}
  ||g||_4^2
  \\
  &\preccurlyeq NT_0L^{\frac{1}{2}}||g||_4^2.
\end{align*}
If $t_g\preccurlyeq \frac{T}{H}$,
then, again, one has $T_0\asymp T\asymp L$ and
$V_0\preccurlyeq\frac{T}{H}$.
By applying \eqref{eqn:triple product} and \eqref{eqn:Omega},
we can bound by $\mathfrak{J}_2$ by \eqref{eqn:reduce to triple product}.
We have deduced that
$\eqref{eqn:reduce to triple product}
\preccurlyeq \frac{T}{H}$
by using Lemma \ref{lem:HK22 lemma 6.1}.
Consequently, we get
\begin{align}\label{eqn:tg small integral -1/2}
  \eqref{eqn:goal Pd -1/2}
  &\preccurlyeq \frac{TH}{C^{3/2}}\left(\frac{N}{T_0L}\right)^{\frac{1}{2}}
  R^{\frac{3}{2}}(t_g+V_0)^{\frac{3}{2}}N^{-\frac{1}{2}}V_0^{-1}
  CV_0\left(\frac{T}{H}\right)^{\frac{1}{2}}\nonumber
  \\
  &\preccurlyeq
  NT^{\frac{3}{2}}\frac{T}{H^{7/2}}
  = o(NT^{\frac{3}{2}}),
\end{align}
which is acceptable.

\subsubsection{$P_{1,\ell}^E(c,r,N)$}
The treatment is very similar to that in \S \ref{subsubsec:contribution d}.
The only difference is that we can not avoid the poles of $B(s,\tau)$ when we shift the path
of integration $\gamma$.
To handle this,
we denote by
\begin{align*}
  \gamma'(\sigma)=\{s=\sigma+iv|\kappa^{3}V_0\leq |v|\leq \kappa^{-3}V_0\}.
\end{align*}
Then we elongate $\gamma(2)$ to $\gamma'(2)$ at the cost of a negligible error.
This allows us to shift $\gamma'(2,V_0)$ to $\gamma'(-\frac{1}{2},V_0)$,
since there is no pole lying on the horizontal line segments.
Then we can do the exact argument as in \S \ref{subsubsec:contribution d},
and see the contribution of the continuous spectrum is acceptable.

In conclusion, we have shown that the contribution of the $\mathcal{J}^+$-term
in \eqref{Voronoi for n} to \eqref{re goal} is bounded by
$NT_0L^{\frac{1}{2}}||g||_4^2+NTH$.
In particular, if $t_g\ll \frac{T^{1+\varepsilon}}{H}$, we actually got that
it can be bounded by
$NT^{\frac{3}{2}}+NTH$, which is not depended on the $L^4$-norm result.

\section{The contribution of $\mathcal{J}^-$}
In this section, we will prove the contribution of $\mathcal{J}^-$
is negligible, and hence complete the proof of Theorem \ref{thm:moment result}.
Note that the $\mathcal{J}^-$-term of \eqref{Voronoi for n} is
\begin{align}\label{eqn: J-}
  \frac{\psi(\bar{d})}{c}\sum_{n\geq 1}\lambda_g(n)e\left(-\frac{\bar{d}n}{c}\right)
  \int_0^\infty V\left(\frac{x}{N}\right)
  e\left(\sigma_1\frac{2\sqrt{mx}}{c}-\sigma_1\frac{t^2c}{4\pi^2\sqrt{mx}}\right)
  \mathcal{J}^-\left(\frac{4\pi\sqrt{nx}}{c}\right)dx,
\end{align}
where
\begin{align*}
  \mathcal{J}^-\left(\frac{4\pi\sqrt{nx}}{c}\right)
  =4\cosh(\pi t_g)K_{2it_g}\left(\frac{4\pi\sqrt{nx}}{c}\right).
\end{align*}
By \eqref{eqn:bound K Bessel}, we can truncate the $n$-sum in \eqref{eqn: J-} at
$n\preccurlyeq\frac{C^2t_g^2}{N}\preccurlyeq N(TH)^{-\frac{1}{2}-\varepsilon}$.
Moreover, by \eqref{eqn:integration rep of K}, we have
\begin{align*}
  4\cosh(\pi t_g)K_{2it_g}(2y)
  &=\sum_{\sigma_2=\pm 1}\int_{-\infty}^\infty \exp(i(2\sigma_2y\sinh v+2t_gv))dv
  \\
  &=\sum_{\sigma_2=\pm 1}\int_{-\infty}^\infty \exp(i(2y\sinh v+2\sigma_2t_gv))dv.
\end{align*}
By partial integration, one may truncate $v$ at $|v|\ll T^\varepsilon$.
Inserting this into \eqref{eqn: J-}, the resulting $x$-integral becomes
\begin{align*}
  \sum_{\sigma_2=\pm 1}\int_0^\infty V\left(\frac{x}{N}\right)
  e\left(\sigma_1\frac{2\sqrt{mx}}{c}-\sigma_1\frac{t^2c}{4\pi^2\sqrt{mx}}
  +\frac{2\sqrt{nx}}{c}\sinh v\right)dx.
\end{align*}
By noting that $\frac{t^2c}{\sqrt{mx}}(\frac{\sqrt{mx}}{c})^{-1}\ll\frac{T^2C^2}{N^2}
\preccurlyeq H^{-2}$, we apply Lemma \ref{lemma:stationary_phase BKY} with
\begin{align*}
  P=U=N, \quad
  X=T^\varepsilon, \quad
  S=\max\left\{\frac{1}{C}, \frac{\sqrt{n}}{N}\frac{|\sinh v|}{C}\right\},\quad
  Y=NS,
\end{align*}
and see that
the $y$-integral is negligible unless $\sgn(v)=\sgn(-\sigma_1)$ and
\begin{align*}
  \frac{\sqrt{nN}}{C}|\sinh v|\asymp \frac{N}{C}.
\end{align*}
So we can smooth $v$ by inserting a nice function $w(v)$ supported on
\begin{align*}
  c_1+\frac{1}{2}\log \frac{N}{n}\leq -\sigma_1 v\leq c_2+\frac{1}{2}\log \frac{N}{n}.
\end{align*}
where $c_1$ and $c_2$ are two constants.
Then the the $v$-integral becomes
\begin{align*}
  \int_{-\infty}^\infty e\left(\frac{2\sqrt{nx}}{c}\sinh v
  +\frac{\sigma_2}{\pi}t_g v\right)w(v)dv.
\end{align*}
It is obvious that, by using $n\preccurlyeq N(TH)^{-\frac{1}{2}-\varepsilon}$,
we have
\begin{align*}
  |\cosh v|\asymp |\sinh v|\asymp \sqrt{\frac{N}{n}}\gg \log T.
\end{align*}
Hence, by using $t_g< T_0\preccurlyeq (TH)^{\frac{3}{4}}$,
$\frac{N}{C}\gg T^{1-\varepsilon}H$ and Lemma \ref{lemma:stationary_phase BKY}
with
\begin{align*}
  X=T^\varepsilon, \quad
  P=U=1, \quad
  Y=S=\frac{N}{C},
\end{align*}
we see that the $v$-integral is negligible
which implies that
the contribution of the $\mathcal{J}^-$ term is negligible.

Therefore, we complete the proof of Theorem \ref{thm:moment result}.


\section{Proof of Corollary \ref{cor:subconvex} and \ref{cor:subconvex2}}
\label{sec:proof of cor}
We only state a sketch proof of Corollary \ref{cor:subconvex2}.
To prove \eqref{eqn:subconvex level=1}, by symmetry, it is natural
to assume $t_g\leq t_f$.
So it follows simply by using \eqref{eqn:main moment 1} with $H=T^{\frac{1}{2}+\varepsilon}$
together with the $L^4$-norm result \eqref{eqn:HK l4}.

Now we prove \eqref{eqn:subconvex burgess level=1}.
Without loss of generality, we assume $q=1$.
By \eqref{eqn:Ki}, we have $||g||_4\ll t_g^\varepsilon$.
If $t_g\ll t_f^{\frac{3}{2}-\varepsilon}$, we use \eqref{eqn:main moment 1}
by letting $H=T_0^{\frac{4}{3}+\varepsilon}T^{-1}$, getting
\eqref{eqn:subconvex burgess level=1}.
If $t_g\gg t_f^{\frac{3}{2}-\varepsilon}$, one has
$t_f\ll t_g^{\frac{2}{3}+\varepsilon}$.
Hence, by \eqref{eqn:main moment 1} with $H=T^{\frac{1}{3}+\varepsilon}$
and switch $f$ and $g$, we deduce that \eqref{eqn:subconvex burgess level=1}
also holds.

This concludes the proof of Corollary \ref{cor:subconvex2}.

\section{Proof of Theorem \ref{thm:shrinkring result for dihedral}}
\label{sec:proof of shrinking result}

We give a sketch proof by closely following the method
in Young \cite[Theorem 1.4 \& Proposition 1.5]{You16}.
Since $t_g$ is the main large parameter now,
we write the abbreviation
\begin{align*}
  A_1\preccurlyeq A_2\Longleftrightarrow A_1\ll_\varepsilon t_g^\varepsilon A_2
\end{align*}
in this section.

Let $h(z)$ be a smooth and compactly supported function
on $\Gamma_0(q)\backslash \mathbb{H}$ which satisfies
\begin{align*}
  ||\Delta^kh||_1\leq C(k)M^{2k},\quad k=0,1,2,\cdots
\end{align*}
where $C(k)$ is a sequence of numbers and $M=M(t_g)$ is a constant
depending on $t_g$.
In particular, the characteristic function of a ball of radius $M^{-1}$
can be approximated by a sequence of $h$.
If $M\leq t_g^\delta$ for some $0<\delta<1$,
then by the Parseval identity and
\begin{align*}
  \left<\phi,h\right>_q\ll \left(\frac{1}{4}+t_\phi^2\right)^{\frac{1}{4}}
  \left(\frac{M^2}{\frac{1}{4}+t_\phi^2}\right)^k
\end{align*}
(see \cite[(4.10)]{You16}), we have
\begin{align}\label{eqn:truncate phi by Young}
  \left<|g|^2,h\right>_q-\left<1,h\right>_q
  =\sum_{\substack{\phi\in \mathcal{B}(q)\\ t_\phi\preccurlyeq M}}
  \left<|g|^2,\phi\right>_q\left<\phi,h\right>_q
  +(Eisenstein\ term) +O(t_g^{-100}).
\end{align}
Using the Cauchy--Schwarz inequality, the first term of the RHS of \eqref{eqn:truncate phi by Young}
can be controlled by
\begin{align*}
  \left(\sum_{\substack{\phi\in \mathcal{B}(q)\\ t_\phi\preccurlyeq M}}
  \left|\left<|g|^2,\phi\right>_q\right|^2\right)^{\frac{1}{2}}||h||_2.
\end{align*}
By Lemma \ref{lem:inner product exchange} and Lemma \ref{lem:Watson-Ichino q},
we have
\begin{align}\label{eqn:triple product in shrinking}
  \sum_{\substack{\phi\in \mathcal{B}(q)\\ t_\phi\preccurlyeq M}}
  \left|\left<|g|^2,\phi\right>_q\right|^2
  \preccurlyeq \sum_{q_1|q}
  \sum_{\substack{\phi\in \mathcal{B}^*(q_1)\\ t_\phi\preccurlyeq M}}
  (t_g(1+|t_\phi|))^{-1}L\left(\frac{1}{2},\phi\right)
  L\left(\frac{1}{2},\phi\times \textup{ad}\,g\right).
\end{align}

Recall $g$ is a dihedral Maass newform.
So we briefly recall some background of
dihedral Maass forms (see \cite[Appendix A.1]{HumKhan20}).
Let $q\equiv 1(\bmod 4)$ be a positive squarefree fundamental discriminant
and let $\psi$ be the primitive quadratic character modulo $q$.
If $g=g_\xi\in \mathcal{B}^*(q,\psi)$ is an dihedral Maass newform associated to
the Hecke Gr$\ddot{\textup{o}}$sscharacter $\xi$ of $\mathbb{Q}(\sqrt{q})$,
then
\begin{align*}
  L\left(\frac{1}{2},\phi\times \textup{ad}\,g_\xi\right)
  =L\left(\frac{1}{2},\phi\times \psi\right)
  L\left(\frac{1}{2},\phi\times g_{\xi^2}\right).
\end{align*}

Now we apply \eqref{eqn:subconvex burgess} to get
$L(\frac{1}{2},\phi\times g_{\xi^2})\ll (t_g+|t_\phi|)^{\frac{3}{4}+\varepsilon}$
and use the standard technique (approximate functional equation $\&$
the large sieve inequality) to obtain
\begin{align*}
  \sum_{\substack{\phi\in \mathcal{B}^*(q_1)\\ t_\phi\preccurlyeq M}}
  L\left(\frac{1}{2},\phi\right)
  L\left(\frac{1}{2},\phi\times \psi\right)\preccurlyeq M^2.
\end{align*}
It follows that
\begin{align*}
  \eqref{eqn:triple product in shrinking}\preccurlyeq
  t_g^{-\frac{1}{4}}M.
\end{align*}
The similar treatments can be done to the Eisenstein term
in \eqref{eqn:truncate phi by Young}.
Actually, we can get a better bound due to Jutila--Motohashi's
uniform subconvexity bound \cite{JutMot05}:
 $L(\frac{1}{2}+it,g)\ll (t_g+|t|)^{1/3+\varepsilon}$.
\footnote{Note that although the subconvexity result
stated in \cite{JutMot05} is for $q=1$, the proof there
can be naturally extended to a general case. In fact, the authors said
``throughout the sequel we shall work with $SL_2(\mathbb{Z})$,
although our argument appears to be effective in a considerably general
setting" (see \cite[p. 62]{JutMot05}).}

Note that $\left<1,h\right>\asymp M^{-2}$ and $||h||_2\asymp M^{-1}$,
so we get \eqref{eqn:shrinking statement} whenever $g$ is a dihedral Maass newform
and $\delta<\frac{1}{12}$.
This completes the proof of Theorem \ref{thm:shrinkring result for dihedral}.

\section*{Acknowledgements}
We would like to thank Bingrong Huang and Yongxiao Lin for helpful comments.



\begin{thebibliography}{10}

\bibitem{Bal}
C.~Balogh,
Asymptotic expansions of the modified Bessel function of the third kind of imaginary order.
\emph{SIAM J. Appl. Math.} 15 (1967), 1315-1323.

\bibitem{BloJanNel}
V.~Blomer, S.~Jana, and P.~Nelson,
The Weyl bound for triple product $L$-functions.
\emph{Duke Math. J.} 172 (2023), no. 6, 1173-1234.

\bibitem{BloBut19}
V.~Blomer and J.~Buttcane,
Global decomposition of $GL(3)$ Kloosterman sums and the spectral large sieve.
\emph{J. Reine Angew. Math.} 757 (2019), 51-88.

\bibitem{BloBut20}
V.~Blomer and J.~Buttcane,
On the subconvexity problem for $L$-functions on $GL(3)$.
\emph{Ann. Sci. $\acute{E}$c. Norm. Sup$\acute{e}$r.} (4) 53 (2020), no. 6, 1441-1500.

\bibitem{BlomerKhanYoung}
V.~Blomer, R.~Khan, and M.~Young,
Distribution of mass of holomorphic cusp forms.
\emph{Duke Math. J.} 162 (2013), no. 14, 2609-2644.

\bibitem{Blomer05}
V.~Blomer,
Rankin-Selberg $L$-functions on the critical line.
\emph{Manuscripta Math.} 117 (2005), no. 2, 111-133.

\bibitem{Bur63}
D. Burgess,
\newblock On character sums and {$L$}-series. {II}.
\newblock \emph{Proc. London Math. Soc.}, (3) 1963, 13:524--536.

\bibitem{DesIwa}
J.M.~Deshouillers, H.~Iwaniec,
Kloosterman sums and Fourier coefficients of cusp forms.
\emph{Invent. Math.} 70, (1982), 219-288.

\bibitem{Dra}
S.~Drappeau,
Sums of Kloosterman sums in arithmetic progressions,
and the error term in the dispersion method.
\emph{Proc. Lond. Math. Soc.} (3) 114 (2017), no. 4, 684–732.

\bibitem{EMOT}
V.~Erd$\acute{\textup{e}}$lyi, W.~Magnus, F.~Oberhettinger, and F.~Tricomi,
Higher transcendental functions I.
\emph{McGraw-Hill.} (1953).

\bibitem{Goo81}
A.~Good,
Cusp forms and eigenfunctions of the Laplacian.
\emph{Math. Ann.} 255 (1981), no. 4, 523-548.

\bibitem{GraRyz}
 I.~S.~Gradshteyn and I.~M.~Ryzhik,
 Tables of Integrals, Series, and Products, 8th ed.,
 \emph{Academic Press}, 2015.

\bibitem{HarMic06}
G.~Harcos and P.~Michel,
The subconvexity problem for Rankin-Selberg $L$-functions and equidistribution
of Heegner points. II.
\emph{Invent. Math.} 163 (2006), no. 3, 581–655.







%




%


\bibitem{HL94} J.~Hoffstein and P.~Lockhart,
Coefficients of Maass forms and the
Siegel zero, with an appendix by Dorian Goldfeld, Hoffstein and Daniel
Lieman.
\emph{Ann. of Math.} (2) 140 (1994), no. 1 , 161-181.

\bibitem{HolSou10}
R.~Holowinsky and K.~Soundararajan
Mass equidistribution for Hecke eigenforms.
\emph{Ann. of Math.} (2) 172 (2010), no. 2, 1517-1528.

\bibitem{HuMicNel}
Y.~Hu, P.~Michel and P.~Nelson,
The subconvexity problem for Rankin-Selberg and triple product $L$-functions.
\emph{ArXiv preprint} (2021), arXiv: 2207.14449.


\bibitem{Hum18}
P.~Humphries,
Equidistribution in shrinking sets and $L^4$-norm bounds for automorphic forms.
\emph{Math. Ann.} 371 (2018), no. 3-4, 1497-1543.

\bibitem{HumKhan22}
P.~Humphries and R.~Khan,
$L^p$-Norm Bounds for Automorphic Forms via Spectral Reciprocity.
\emph{Proc. Lond. Math. Soc.} (3) 130 (2025), no. 6, Paper No. e70061, 80 pp.

\bibitem{HumKhan20}
P.~Humphries and R.~Khan,
On the random wave conjecture for dihedral Maass forms.
\emph{Geom. Funct. Anal.} 30 (2020), no. 1, 34–125.

\bibitem{Ich08}
A.~Ichino,
Trilinear forms and the central values of triple product $L$-functions.
\emph{Duke Math. J.} (145) 2 (2008), 281-307.

\bibitem{Iwa90} H.~Iwaniec,
Small eigenvalues of Laplacian for $\Gamma_0(N)$.
\emph{Acta Arith.} 56 (1990),
no. 1, 65–82.











\bibitem{IwaKow04}
H.~Iwaniec and E.~Kowalski,
Analytic number theory, volume~53 of
\emph{American Mathematical
  Society Colloquium Publications}.
American Mathematical Society, Providence, RI, 2004.

\bibitem{IwaLuoSar00}
H.~Iwaniec; W.~Luo and P.~Sarnak,
Low lying zeros of families of $L$-functions.
\emph{Inst. Hautes Études Sci. Publ. Math.}
(91) 1 (2000), 55-131.

\bibitem{Jut96}
M.~Jutila,
The additive divisor problem and its analogs for Fourier coefficients of cusp forms. I.
\emph{Math. Z.} 223 (1996), no. 3, 435-461.

\bibitem{Jut97}
M.~Jutila,
The additive divisor problem and its analogs for Fourier coefficients of cusp forms. II.
\emph{Math. Z.} 225 (1997), no. 4, 625-637.

\bibitem{JutMot05}
 M.~Jutila and Y.~Motohashi,
Uniform bound for Hecke $L$-functions.
\emph{Acta Math.} 195 (2005), 61-115.

\bibitem{JutMot06}
 M.~Jutila and Y.~Motohashi,
 Uniform bounds for Rankin-Selberg $L$-functions.
 \emph{Multiple Dirichlet series, automorphic forms, and analytic number theory}.
  243–256, Proc. Sympos. Pure Math., 75, Amer. Math. Soc., Providence, RI, 2006.

\bibitem{Ki}
 H.~Ki, $L^4$-norms and sign changes of Maass forms,
 \emph{ArXiv preprint} (2023), arXiv: 2302.02625.


\bibitem{Kim2003}
H.~Kim,
Functoriality for the exterior square of $\GL_4$ and the symmetric fourth of $\GL_2$,
\emph{J. Amer. Math. Soc.} 16 (2003), no. 1, 139--183.
With appendix 1 by Ramakrishnan and appendix 2 by Kim and Sarnak.

\bibitem{KirPetYou19}
E.~K{\i}ral, I.~Petrow, and M.~Young,
Oscillatory integrals with uniformity in parameters.
\emph{J. Th\'eor. Nombres Bordeaux} 31 (2019), no. 1, 145--159.

\bibitem{KowMicVan02}
E.~Kowalski; P.~Michel and J.~VanderKam,
Rankin-Selberg $L$-functions in the level aspect.
\emph{Duke Math. J.} 114 (2002), no. 1, 123-191.

\bibitem{KroSta04}
B.~Kr$\ddot{\textup{o}}$tz and Bernhard; R.~Stanton,
Holomorphic extensions of representations. I. Automorphic functions.
\emph{Ann. of Math.} (2) 159 (2004), no. 2, 641-724.








\bibitem{LauLiuYe06}
Y-K.~Lau; J.~Liu and Y.~Ye,
 A new bound $k^{2/3+\varepsilon}$
 for Rankin-Selberg L-functions for Hecke congruence subgroups.
 \emph{IMRP Int. Math. Res. Pap.} 2006, Art. ID 35090, 78 pp.

\bibitem{Li10}
X.~Li,
Upper bounds on L-functions at the edge of the critical strip.
\emph{Int. Math. Res. Not. IMRN} 2010, no. 4, 727-755.

\bibitem{Li11}
X.~Li,
\newblock Bounds for {${\rm GL}(3)\times {\rm GL}(2)$} {$L$}-functions and
  {${\rm GL}(3)$} {$L$}-functions.
\newblock {\em Ann. of Math. (2)}, 173(1):301--336, 2011.

\bibitem{LiYou12}
X.~Li and M.~Young,
The $L^2$ restriction norm of a $GL_3$ Maass form.
\emph{Compos. Math.} 148 (2012), no. 3, 675-717.

\bibitem{Lin06}
E.~Lindenstrauss,
Invariant measures and arithmetic quantum unique ergodicity.
\emph{Ann. of Math.} (2) 163 (2006), no. 1, 165-219.


\bibitem{Luo}
W.~Luo,
$L^4$-norms of the dihedral Maass forms.
\emph{Int. Math. Res. Not. IMRN} (2014), no. 8, 2294-2304.









%


%
%


%
\bibitem{Mic04}
P.~Michel,
The subconvexity problem for Rankin-Selberg $L$-functions
and equidistribution of Heegner points.
\emph{Ann. of Math.} (2) 160 (2004), no. 1, 185-236.



\bibitem{MicVen10}
P.~Michel and A.~Venkatesh,
The subconvexity problem for $\GL_2$.
\emph{Publ. Math. Inst. Hautes \'Etudes Sci.} No. 111 (2010), 171--271.




%



\bibitem{Mun15III}
R.~Munshi,
The circle method and bounds for $L$-functions--III: $t$-aspect subconvexity for $\GL(3)$ $L$-functions.
\emph{J. Amer. Math. Soc.} 28 (2015), no. 4, 913--938.


\bibitem{Mun15IV}
R.~Munshi,
The circle method and bounds for $L$-functions--IV: Subconvexity for twists of $\GL(3)$ $L$-functions. \emph{Ann. of Math. (2)} 182 (2015), no. 2, 617--672.

\bibitem{Nel21}
P.~Nelson,
Bounds for standard $L$-functions.
\emph{ArXiv preprint} (2021), arXiv:2109.15230.





\bibitem{PetYou20}
I. Petrow and M. Young,
The Weyl bound for Dirichlet $L$-functions of cube-free conductor.
\emph{Ann. of Math. (2)} 192 (2020), no. 2, 437--486. 

\bibitem{PetrowYoung2019}
I. Petrow and M. Young,
The fourth moment of Dirichlet $L$-functions along a coset and the Weyl bound.
\emph{ArXiv preprint} (2019), arXiv:1908.10346.
to appear in  Duke Math. J..

\bibitem{RudSar94}
Z.~Rudnick and P.~Sarnak,
The behaviour of eigenstates of arithmetic hyperbolic manifolds.
\emph{Comm. Math. Phys.} 161 (1994), no. 1, 195-213.

\bibitem{Sar01}
P.~Sarnak,
Estimates for Rankin-Selberg $L$-functions and quantum unique ergodicity.
\emph{J. Funct. Anal.} 184 (2001), no. 2, 419-453.

\bibitem{Sou10}
K.~Soundararajan,
Quantum unique ergodicity for $\SL_2(\mathbb{Z})\backslash \mathbb{H}$.
\emph{Ann. of Math.} (2) 172 (2010), no. 2, 1529-1538.


\bibitem{Tit86}
E.~C.~Titchmarsh, The Theory of the Riemann Zeta-Function, 2nd ed.,
\emph{Clarendon Press,}
Oxford Univ. Press, New York, 1986.


\bibitem{Wat02}
T.~Watson,
Rankin triple products and quantum chaos, Ph.D. Thesis,
\emph{Princeton University,} 2002 (revised 2008).


\bibitem{Weyl21}
H. Weyl,
Zur absch\"atzung von $\zeta(1+it)$.
\emph{Math. Z.}, 10 (1921), 88--101.

\bibitem{You16}
M.~Young,
The quantum unique ergodicity conjecture for thin sets.
\emph{Adv. Math.} 286 (2016), 958-1016.

\bibitem{Zac}
R. Zacharias,
Simultaneous non-vanishing for Dirichlet L-functions.
\emph{Ann. Inst. Fourier (Grenoble)} 69 (2019),
no. 4, 1459–1524.

\end{thebibliography}
\end{document}